\documentclass{amsart}
\usepackage{amsmath,amssymb,amsthm}
\usepackage{hyperref}

\numberwithin{equation}{section}
\newtheorem{theorem}{Theorem}[section]
\newtheorem{lemma}[theorem]{Lemma}
\newtheorem{proposition}[theorem]{Proposition}
\newtheorem{definition}[theorem]{Definition}
\newtheorem{corollary}[theorem]{Corollary}

\newcommand{\orcid}[1]{\href{https://orcid.org/#1}{\texttt{ORCID: #1}}}

\begin{document}

\title[Small caps, exponential sums and lattice points]{%
Small cap decoupling, exponential sums\\
and lattice points}

\date{}

\author{Xiangyu Wang} \address{Xiangyu Wang \\ Department of Mathematics, University of Illinois Urbana-Champaign\\ \newline \orcid{0009-0003-5983-6961}} \email{xw70@illinois.edu}

\begin{abstract}

We apply the recent small cap decoupling inequalities of Jeong and Oh for the
moment curve to exponential sums. For sums in one variable we prove a new
estimate under conditions on two consecutive derivatives, and we deduce a
family of new exponent pairs. For the Gauss circle
problem and the Dirichlet divisor problem we improve the exponent in the
bounds for the remainder terms to $293/940=0.3117\ldots$.

\end{abstract}


\maketitle

\section{Introduction}

Let $d\ge 2$ and $J\ge 1$ be integers. For complex numbers $b_1,\dots,b_J$ we consider the Weyl sums
\begin{equation}\label{eq:W}
W_b(y)=\sum_{r=1}^{J}b_r\,e\bigl(y_1r+y_2r^2+\dots+y_dr^d\bigr)\qquad(y\in\mathbb{R}^d),
\end{equation}
where $e(z)=\exp(2\pi iz)$. The mean values of $|W_b|$ over the torus $\mathbb{T}^d$ are described by Vinogradov's mean value theorem; see \cite{BDG}. Small cap decoupling, introduced by Demeter, Guth and Wang~\cite{DGW}, concerns thinner sets, such as $\mathbb{T}^{d-1}\times I$ with a short interval $I$. The sums \eqref{eq:W} correspond to the moment curve $(t,t^2,\dots,t^d)$, for which sharp small cap decoupling inequalities were proved by Guth and Maldague~\cite{GM} in $\mathbb{R}^3$, and recently by Jeong and Oh~\cite{JO} in $\mathbb{R}^d$ for every $d\ge 4$.

In this paper we apply the results of \cite{JO} to two problems on exponential sums. The mechanism is the same in both. When the phase is expanded in a Taylor polynomial of degree $d$ around short shifts of the variable of summation, the sum becomes a bilinear form in the Weyl sums \eqref{eq:W}, evaluated at the vectors of Taylor coefficients. If the phase is smooth, the coefficient of degree $d$ is small. The vectors of coefficients then lie in a set of the form $\mathbb{T}^{d-1}\times I$ with a short interval $I$, and on such sets small cap decoupling gives more than Vinogradov's mean value theorem.

\subsection{Exponential sums in one variable}

Our first result is an estimate for sums $\sum e(f(n))$ in which two consecutive derivatives of $f$ have a given size.

\begin{theorem}\label{thm:onevar}
Let $d\ge 4$ be an integer, and put $\mu_d=d(d-1)$ if $d\ge 5$ and $\mu_4=14$. Let $\mathcal{K}$ be a compact subinterval of $(\frac{d+1}{2},\,d-\frac1d)$, and let $0<c\le C$. Let $N\ge 2$ and $t\in\mathcal{K}$, let $I\subseteq[N,2N]$ be an interval, and let $f:I\to\mathbb{R}$ be a function with $d+1$ continuous derivatives such that
\[
cN^{t-i}\le|f^{(i)}(x)|\le CN^{t-i}\qquad(x\in I,\ i\in\{d,d+1\}).
\]
Then
\begin{equation}\label{eq:onevar}
\sum_{n\in I\cap\mathbb{Z}}e\bigl(f(n)\bigr)\lesssim_{d,\mathcal{K},c,C,\varepsilon}N^{1+\varepsilon}\Bigl(N^{-\frac{1}{\mu_d}}+N^{-\frac{d+1-t}{d(d+1)}}\Bigr).
\end{equation}
\end{theorem}

The numbers $N^{t-d}$ and $N^{t-d-1}$ are the sizes of $f^{(d)}$ and $f^{(d+1)}$. For comparison, let $t$ be as in the theorem. Then the estimate of Heath-Brown \cite[Theorem~1]{HB} for the derivative of order $d$ gives the bound
\[
N^{1+\varepsilon}\Bigl(N^{-\frac{d-t}{d(d-1)}}+N^{-\frac{1}{d(d-1)}}\Bigr),
\]
and for the derivative of order $d+1$ it gives
\[
N^{1+\varepsilon}\Bigl(N^{-\frac{d+1-t}{d(d+1)}}+N^{-\frac{1}{d(d+1)}}\Bigr).
\]
If $d\ge 5$, the right-hand side of \eqref{eq:onevar} consists of the second term of the first bound and the first term of the second bound. In this case \eqref{eq:onevar} agrees with the first bound for $t\le d-1$, and it is smaller than both bounds for $d-1<t<d-\frac1d$. For example, if $d=5$ and $t=\frac92$, the two bounds of Heath-Brown are $N^{39/40+\varepsilon}$ and $N^{29/30+\varepsilon}$, while \eqref{eq:onevar} gives $N^{19/20+\varepsilon}$. If $d=4$, the estimate \eqref{eq:onevar} is smaller than both bounds for $\frac{22}{7}<t<\frac{15}{4}$.

Theorem~\ref{thm:onevar} gives new exponent pairs. We use this notion in the sense of the Analytic Number Theory Exponent Database \cite[Definition~5.1]{DB}; the definition is recalled in Section~\ref{subsec:pairs}. Heath-Brown \cite[Theorem~2]{HB} proved that
\begin{equation}\label{eq:hbpairs}
(p_m,q_m)=\Bigl(\frac{2}{(m-1)^2(m+2)},\ 1-\frac{3m-2}{m(m-1)(m+2)}\Bigr)
\end{equation}
is an exponent pair for every integer $m\ge 3$; see also \cite[Theorem~5.17]{DB}.

\begin{corollary}\label{cor:pairs}
For every integer $d\ge 4$,
\begin{equation}\label{eq:kl}
(k_d,\ell_d)=\Bigl(\frac{2}{(d+2)(d^2+d+1)},\ 1-\frac{3d^2+3d-1}{(d+1)(d+2)(d^2+d+1)}\Bigr)
\end{equation}
is an exponent pair. Moreover, $k_d<p_{d+1}$ and $\ell_d<q_{d+1}$.
\end{corollary}

For example, $(k_4,\ell_4)=(\frac{1}{63},\frac{571}{630})$ and $(k_5,\ell_5)=(\frac{2}{217},\frac{1213}{1302})$, while $(p_5,q_5)=(\frac{1}{56},\frac{127}{140})$ and $(p_6,q_6)=(\frac{1}{100},\frac{14}{15})$. Thus Corollary~\ref{cor:pairs} improves the pairs $(p_m,q_m)$ with $m\ge 5$ in both coordinates. One more new exponent pair, which is not of the form \eqref{eq:kl}, is given in Proposition~\ref{prop:star}. In Section~\ref{subsec:hull} we compare the new pairs with the exponent pairs that are recorded in \cite[Chapter~5]{DB}. All the new pairs lie outside the convex hull of the recorded ones. The pairs $(p_m,q_m)$ with $m\ge 6$ are vertices of that convex hull, but they are no longer extreme points when the new pairs are added.

\subsection{The circle and divisor problems}

For real $X\ge 4$ let
\[
\begin{aligned}
E(X)&=\#\{(u,v)\in\mathbb{Z}^2:\ u^2+v^2\le X\}-\pi X,\\
\Delta(X)&=\#\{(a,b)\in\mathbb{N}^2:\ ab\le X\}-X\log X-(2\gamma-1)X,
\end{aligned}
\]
where $\mathbb{N}$ is the set of positive integers and $\gamma$ is Euler's constant. These are the remainder terms in the Gauss circle problem and in the Dirichlet divisor problem; the number of pairs $(a,b)$ in the second line is the sum of the numbers of divisors of the integers $n\le X$. It is conjectured that both remainder terms are $O_\varepsilon(X^{1/4+\varepsilon})$.

The exponent $1/3$ goes back to Vorono\"{\i} and Sierpi\'nski, and it was lowered many times by van der Corput's method. The strongest results come from the method of Bombieri and Iwaniec~\cite{BI}, which Iwaniec and Mozzochi~\cite{IM} adapted to these two problems. Huxley~\cite{Huxley} reached the exponent $131/416=0.31490\ldots$, and Li and Yang~\cite{LY} improved it to $0.31448\ldots$; both results hold for $E(X)$ and for $\Delta(X)$. Recently Cen, Gao, Zhang and Zhang~\cite{CGZZ} proved that $E(X)\lesssim_\varepsilon X^{209/670+\varepsilon}$, where $209/670=0.311940\ldots$. We combine their argument with the results of~\cite{JO}.

\begin{theorem}\label{thm:main}
For every $\varepsilon>0$ and every real $X\ge 4$,
\[
|E(X)|+|\Delta(X)|\lesssim_\varepsilon X^{293/940+\varepsilon}.
\]
\end{theorem}

Here $293/940=0.311702\ldots$. The divisor problem is not treated in \cite{CGZZ}, but the main estimate of that paper also gives $\Delta(X)\lesssim_\varepsilon X^{209/670+\varepsilon}$, by the argument of Section~\ref{subsec:proofmain}.

\subsection{Outline of the proofs}

The moments of the sums \eqref{eq:W} over $\mathbb{T}^{d-1}\times I$ are bounded in \cite[Theorem~1.3]{JO}. For our purposes it is more efficient to work with level sets. From the proof of the small cap decoupling theorem in \cite{JO} we extract the following estimate (Proposition~\ref{prop:levelset}): if $d\ge 4$, $|b_r|\le 1$ and $J^{-(d-1)}\le|I|\le 1$, then for every $\lambda>0$
\[
\bigl|\{y\in\mathbb{T}^{d-1}\times I:\ |W_b(y)|>\lambda\}\bigr|\lesssim_{d,\varepsilon}J^{\varepsilon-d(d+1)/2}\Bigl(|I|\Bigl(\frac J\lambda\Bigr)^{d(d+1)}+\Bigl(\frac J\lambda\Bigr)^{\mu_d}\Bigr),
\]
with $\mu_d$ as in Theorem~\ref{thm:onevar}. The first term corresponds to square root cancellation, and for $d\ge 5$ the second term corresponds to constructive interference. In Section~\ref{sec:weyl} we combine this estimate with the double large sieve of Bombieri and Iwaniec~\cite{BI}. The result is Lemma~\ref{lem:ls}, which bounds a bilinear form in the values of $W_b$ at given points in terms of the number of pairs of points that are close to each other. Together with Lemma~\ref{lem:shift}, which describes the passage from exponential sums to Weyl sums by short shifts, it is the common basis of the two applications.

In Theorem~\ref{thm:onevar} the Taylor coefficient of degree $d$ of $f$ has size $N^{t-d}$, and Lemma~\ref{lem:ls} is applied with an interval $I$ of this length. The pairs of close coefficient vectors are controlled by the two coefficients of the highest degrees: both are monotonic in $n$, and for $t$ in the range of the theorem they separate any two different integers $n$ (Lemma~\ref{lem:sep}). This is why two derivatives occur in Theorem~\ref{thm:onevar}. The argument is a variant of Vinogradov's method in the form given by Heath-Brown~\cite{HB}, with Proposition~\ref{prop:levelset} in place of Vinogradov's mean value theorem.

For the circle and divisor problems, both remainder terms are written as sums of the sawtooth function $\psi(x)=x-\lfloor x\rfloor-\frac12$, and these sums are reduced to double exponential sums by Fourier expansion. This reduction is classical. We recall it in Section~\ref{subsec:reduction}, where we also introduce the notation of \cite{CGZZ} and state the main estimate for the exponential sums, Theorem~\ref{thm:rec}. The proof of Theorem~\ref{thm:rec} follows \cite{CGZZ}. In the critical range of the parameters, the exponential sums are first transformed by the van der Corput $B$-process. The new variable of summation is then shifted, and the phase is expanded up to degree six. In \cite{CGZZ} the vectors of the six Taylor coefficients lie in the product of $\mathbb{T}^5$ and a short interval; the sixth coordinate is held fixed, and the size of the Weyl sums is controlled by Vinogradov's mean value theorem in degree five. Our starting point is that the fifth coefficient also varies in a short interval, whose length is a negative power of $X$. Thus Lemma~\ref{lem:ls} can be applied with $d=5$; the lemma allows one more coordinate, which is used for the sixth coefficient. The main term of the resulting bound is smaller than in \cite{CGZZ} by a negative power of $X$ (see Lemma~\ref{lem:csieve}), so that the large sieve becomes effective in a larger range of the parameters. In this application the count of the close pairs is the difficult part of \cite{CGZZ}. This count and the estimate of \cite{CGZZ} for the remaining range are proved there for other values of the parameters. We list the properties of the parameters that their proofs use, and we verify them for our values (Tables~\ref{tab:close} and~\ref{tab:watt}). In Section~\ref{subsec:sharp} we show that, for this argument, the choice of the shift length is optimal and the exponent $293/940$ cannot be lowered by changing the parameters.

\subsection*{Structure of the paper}
Section~\ref{sec:weyl} contains the level set estimate, the large sieve for localized Weyl sums and the lemma on short shifts. Theorem~\ref{thm:onevar} and Corollary~\ref{cor:pairs} are proved in Section~\ref{sec:onevar}, where we also compare the new exponent pairs with the known ones. Theorem~\ref{thm:main} is proved in Section~\ref{sec:lattice}.

\subsection*{Notation}
We write $A\lesssim B$ if $A\le CB$ for a constant $C$; a subscript indicates the permitted dependence of $C$. We write $A\asymp B$ if $A\lesssim B\lesssim A$, and $n\approx M$ if $M\le n<2M$. Every estimate that involves $\varepsilon$ holds for each $\varepsilon>0$. Further, $e(z)=\exp(2\pi iz)$, $\|x\|=\min_{j\in\mathbb{Z}}|x-j|$ for real $x$, and $\mathbb{T}=\mathbb{R}/\mathbb{Z}$. The sums $W_b$ are defined by \eqref{eq:W}; their degree $d$ and their length $J$ will be clear from the context. In Section~\ref{sec:lattice} we use the notation of \cite{CGZZ}. A few letters therefore have a different meaning there than in Section~\ref{sec:onevar}; this is pointed out where it occurs.

\subsection*{AI Declaration}
In an earlier draft, the author obtained a smaller improvement over \cite{CGZZ} in the circle and divisor problems from the moment bounds of \cite[Theorem~1.3]{JO}. Following the author's instructions, GPT-6.2 Astra extracted from this argument a separate note on exponential sums in one variable, with estimates that are weaker than Theorem~\ref{thm:onevar} and with some exponent pairs that are superseded by those given here. Claude Fable 5.2 strengthened these arguments. It derived the level set estimate of Proposition~\ref{prop:levelset} from \cite[Lemmas~12.3 and~12.4]{JO} and used it in place of the moment bounds; this gives the exponent $293/940$ in Theorem~\ref{thm:main} and Theorem~\ref{thm:onevar} in its present form. It also worked out the exponent pairs of Corollary~\ref{cor:pairs} and Proposition~\ref{prop:star}. In addition, Claude Fable 5.2 assisted with the computations, with the verification of the proofs, and with language editing. The example in Section~\ref{subsec:sharp} was constructed by GPT-6.2 Astra at the author's request, after the author had tried several choices of the shift length $J$ without improving the final exponent.

\section{Localized Weyl sums}\label{sec:weyl}

In Sections~\ref{subsec:levelsets} and~\ref{subsec:ls}, $d\ge 4$ and $J\ge 1$ are integers, and $W_b$ is the Weyl sum \eqref{eq:W} of degree $d$ and length $J$. As in Theorem~\ref{thm:onevar}, we put $\mu_d=d(d-1)$ if $d\ge 5$ and $\mu_4=14$.

\subsection{Level sets}\label{subsec:levelsets}

After rescaling, $W_b$ becomes a function on $\mathbb{R}^d$ whose Fourier transform is concentrated near the points $(t,t^2,\dots,t^d)$ of the moment curve, and this is the setting of small cap decoupling. The following estimate is not stated in \cite{JO}, but it follows from the proof of the small cap decoupling theorem in that paper.

\begin{proposition}\label{prop:levelset}
Let $I\subset\mathbb{R}$ be an interval with $J^{-(d-1)}\le|I|\le 1$. Then, for all complex numbers $b_1,\dots,b_J$ with $|b_r|\le 1$ and all $\lambda>0$,
\begin{equation}\label{eq:levelset}
\bigl|\{y\in\mathbb{T}^{d-1}\times I:\ |W_b(y)|>\lambda\}\bigr|\lesssim_{d,\varepsilon}J^{\varepsilon-d(d+1)/2}\Bigl(|I|\Bigl(\frac J\lambda\Bigr)^{d(d+1)}+\Bigl(\frac J\lambda\Bigr)^{\mu_d}\Bigr).
\end{equation}
\end{proposition}

For $d=5$ the right-hand side is $J^{\varepsilon}(|I|J^{15}\lambda^{-30}+J^{5}\lambda^{-20})$. If $d\ge 5$ and $|I|=J^{-\sigma}$, integration of \eqref{eq:levelset} over $\lambda$ gives back the moment bounds of \cite[Theorem~1.3]{JO} for the exponents between $d(d-1)$ and $d(d+1)$.

\begin{proof}
We use the notation of \cite[Section~12]{JO} for the moment curve in $\mathbb{R}^d$, except that the level is denoted by $\lambda$. To shorten the formulas we write $u=d(d-1)/2$ and $v=d(d+1)/2$ in this proof, so that $v-u=d$. Write $|I|=J^{-\sigma}$ with $0\le\sigma\le d-1$, and let $R=J^{d-\sigma}$. The small caps $\gamma$ are those of width $J^{-1}$. In \cite{JO} this width is written as a power of $R$; the exponent is $1/(d-\sigma)$ here, and the implied constants in \cite{JO} may depend on it. This is harmless: we may enlarge $I$ so that $\sigma$ becomes an integer multiple of $\varepsilon$, which costs at most a factor $J^{\varepsilon}$ in \eqref{eq:levelset}, and then only finitely many values of $\sigma$ occur. We may also assume that $\lambda\ge 2J^{1/2}$ and that $b_J=0$. Indeed, for $\lambda\le 2J^{1/2}$ the estimate \eqref{eq:levelset} is trivial, and for larger $\lambda$ the term with $r=J$ can be removed at the cost of replacing $\lambda$ by $\lambda/2$.

\emph{Reduction to small caps.} This step is the same as in the deduction of \cite[Theorem~1.3]{JO} from \cite[Theorem~1.4]{JO}; see \cite[Section~1.1]{GM}. Let $\varphi$ be a Schwartz function on $\mathbb{R}^d$ with $|\varphi|\ge 1$ on $[0,1]^d$ whose Fourier transform is supported in $[0,c]^d$, where $c>0$ is a small constant that depends only on $d$. Let $x_d^0=J^d\inf I$, and for $1\le r<J$ put
\[
f_r(x)=b_r\,e\Bigl(\sum_{i=1}^{d}x_i\Bigl(\frac rJ\Bigr)^{i}\Bigr)\,\varphi\Bigl(\frac{x_1}{J^d},\dots,\frac{x_{d-1}}{J^d},\frac{x_d-x_d^0}{R}\Bigr)\qquad(x\in\mathbb{R}^d).
\]
The Fourier transform of $f_r$ is supported in the translate of $[0,cJ^{-d}]^{d-1}\times[0,cR^{-1}]$ by the point of the moment curve with parameter $r/J$. Since $R\le J^d$, this set is contained in one small cap if $c$ is small enough, and different values of $r$ give different caps. Moreover, $\|f_r\|_\infty\le\|\varphi\|_\infty$ and $\|f_r\|_2^2\lesssim J^{d(d-1)}R$. Let $f=\sum_rf_r$. If $x\in[0,J^d)^{d-1}\times J^dI$ and $y_i=x_i/J^i$, then $|f(x)|\ge|W_b(y)|$. The function $W_b$ has period $1$ in each of $y_1,\dots,y_{d-1}$. Hence each point of $\mathbb{T}^{d-1}\times I$ corresponds to $J^{u}$ points $x$, and $dx=J^{v}\,dy$. Since $u+v=d^2$, it follows that the left-hand side of \eqref{eq:levelset} is at most
\[
J^{-d^2}\,\bigl|\{x\in\mathbb{R}^d:\ |f(x)|>\lambda\}\bigr| .
\]

\emph{The level set of $f$.} We claim that, if $f=\sum_\gamma f_\gamma$, where the Fourier transform of $f_\gamma$ is supported in $\gamma$ and $\|f_\gamma\|_\infty\le 1$, then
\begin{equation}\label{eq:levelsetf}
\bigl|\{|f|>\lambda\}\bigr|\lesssim_{d,\varepsilon}R^{\varepsilon}\Bigl(J^{v-1}\lambda^{-2v}+J^{u-1}R^{-1}\lambda^{-2u}\Bigl(\frac J\lambda\Bigr)^{\mu_d-2u}\Bigr)\sum_\gamma\|f_\gamma\|_2^2 .
\end{equation}
Note that $\mu_d-2u$ is equal to $0$ if $d\ge 5$, and to $2$ if $d=4$. We apply \eqref{eq:levelsetf} to the function $f/\|\varphi\|_\infty$, with $f$ as constructed above. For it the sum over $\gamma$ is $O(J^{2u+1}R)$. Since $R=J^{d}|I|$ and $u+v=d^2$, the proposition follows from \eqref{eq:levelsetf}.

To prove \eqref{eq:levelsetf} we group the small caps into the canonical caps $\theta\in\Theta(R^{-1/d})$ and write $f_\theta=\sum_{\gamma\subset\theta}f_\gamma$. As in \cite[Lemma~12.4]{JO}, the letter $\nu$ denotes a cap of width $|\nu|$, where $|\nu|$ is a dyadic number with $R^{-1/d}\le|\nu|\le 1$ and $\nu$ belongs to the family $\Theta(|\nu|)$. For such a cap we put $g_\nu=\sum_{\theta\subset\nu}|f_\theta|^2$. By \cite[Lemma~12.3]{JO}, which is a consequence of the amplitude dependent wave envelope estimate \cite[Theorem~1.7]{JO},
\[
\lambda^{v+1}\,\bigl|\{|f|>\lambda\}\bigr|\le C_\varepsilon R^{\varepsilon}\sum_{\nu}\int_{\{g_\nu>C_\varepsilon^{-1}R^{-\varepsilon}\lambda^2|\nu|^2\}}g_\nu^{(v+1)/2},
\]
where the sum is over all these caps. We have $u\ge(v+1)/2$, and on the domain of integration $g_\nu^{(v+1)/2}\le(C_\varepsilon R^{\varepsilon}\lambda^{-2}|\nu|^{-2})^{u-(v+1)/2}g_\nu^{u}$. Hence
\begin{equation}\label{eq:after123}
\lambda^{2u}\,\bigl|\{|f|>\lambda\}\bigr|\lesssim_{d,\varepsilon}R^{d^2\varepsilon}\sum_{\nu}|\nu|^{-(2u-v-1)}\int_{\{g_\nu>C_\varepsilon^{-1}R^{-\varepsilon}\lambda^2|\nu|^2\}}g_\nu^{u}.
\end{equation}
By \cite[Lemma~12.4]{JO}, applied to $\nu$ and to the caps $\theta\subset\nu$ as in \cite[Section~12.2]{JO}, there is a constant $C$ such that
\[
\begin{aligned}
\int_{\mathbb{R}^d}g_\nu^{u}&\lesssim_{d,\varepsilon}R^{\varepsilon}(J|\nu|)^{u-1}\sum_{\gamma\subset\nu}\|f_\gamma\|_2^2,\\
\int_{\{g_\nu>CJ|\nu|\}}g_\nu^{u}&\lesssim_{d,\varepsilon}R^{\varepsilon}(J|\nu|)^{u-1}\Bigl(\frac{R^{-1/d}}{|\nu|}\Bigr)^{u/2}\sum_{\gamma\subset\nu}\|f_\gamma\|_2^2 .
\end{aligned}
\]
For the caps with $C_\varepsilon^{-1}R^{-\varepsilon}\lambda^2|\nu|^2\ge CJ|\nu|$ we use the second bound. Together with the power of $|\nu|$ in \eqref{eq:after123} it gives the factor
\[
|\nu|^{-(2u-v-1)}(J|\nu|)^{u-1}\Bigl(\frac{R^{-1/d}}{|\nu|}\Bigr)^{u/2}=J^{u-1}R^{-\frac{d-1}{4}}\,|\nu|^{-\frac{d(d-5)}{4}}.
\]
If $d\ge 5$, the exponent of $|\nu|$ is not positive. Since $|\nu|\ge R^{-1/d}$, the factor is then at most $J^{u-1}R^{-\frac{d-1}{4}+\frac{d-5}{4}}=J^{u-1}R^{-1}$. If $d=4$, the factor is $J^{5}R^{-3/4}|\nu|$. In this case we also use the pointwise bound $g_\nu\lesssim J^{2}R^{-1/4}|\nu|$ of \cite[Lemma~12.4(i)]{JO}. It shows that the domain of integration in \eqref{eq:after123} is empty unless $|\nu|\lesssim C_\varepsilon R^{\varepsilon}J^{2}R^{-1/4}\lambda^{-2}$, and for these caps the factor is $O(C_\varepsilon R^{\varepsilon}J^{7}R^{-1}\lambda^{-2})$; here $J^{7}R^{-1}\lambda^{-2}=J^{u-1}R^{-1}(J/\lambda)^{2}$. There are $O(\log R)$ dyadic widths, so in both cases the contribution of these caps to the right-hand side of \eqref{eq:after123} is $O(R^{2d^2\varepsilon}J^{u-1}R^{-1}(J/\lambda)^{\mu_d-2u}\sum_\gamma\|f_\gamma\|_2^2)$. The other caps have $|\nu|<CC_\varepsilon R^{\varepsilon}J\lambda^{-2}$. For them we use the first bound, and since $|\nu|^{-(2u-v-1)}(J|\nu|)^{u-1}=J^{u-1}|\nu|^{d}$, their contribution is $O(R^{2d^2\varepsilon}J^{v-1}\lambda^{-2d}\sum_\gamma\|f_\gamma\|_2^2)$. Since $2u+2d=2v$, this proves \eqref{eq:levelsetf}.
\end{proof}

\subsection{A large sieve for localized Weyl sums}\label{subsec:ls}

The double large sieve of Bombieri and Iwaniec~\cite{BI} bounds a bilinear form in the values of a Weyl sum at given points in terms of the number of pairs of points that are close to each other and of the size of the Weyl sum. We need it in the following form, in which the size enters through Proposition~\ref{prop:levelset}. The lemma is stated for sums of degree $d+1$, because this is what Section~\ref{sec:lattice} requires; the case of the sums $W_b$ of degree $d$ is included. The proof follows that of \cite[Lemma~B.1]{CGZZ}; see also \cite{BW}.

\begin{lemma}\label{lem:ls}
Let $D$ and $D'$ be real numbers with $J^{-(d-1)}\le D\le\frac14$ and $D'\ge 0$. Let $x_1,\dots,x_L$ be points of $\mathbb{R}^{d+1}$ with $|x_{j,d}|\le D$ and $|x_{j,d+1}|\le D'$. Let $B$ be the number of pairs $(j,j')\in\{1,\dots,L\}^2$ such that
\[
\|x_{j,i}-x_{j',i}\|\le J^{-i}\quad(1\le i<d),\qquad|x_{j,i}-x_{j',i}|\le J^{-i}\quad(i\in\{d,d+1\}),
\]
and put $\Delta=1+D'J^{d+1}$. Then, for all complex numbers $a_j$ and $b_r$ with $|a_j|\le 1$ and $|b_r|\le 1$,
\begin{equation}\label{eq:ls}
\begin{aligned}
\Bigl|\sum_{j=1}^{L}a_j\sum_{r=1}^{J}b_r\,e\Bigl(\sum_{i=1}^{d+1}x_{j,i}\,r^i\Bigr)\Bigr|\lesssim_{d,\varepsilon}J^{1+\varepsilon}\Bigl[&(\Delta D)^{\frac{1}{d(d+1)}}L^{1-\frac{2}{d(d+1)}}B^{\frac{1}{d(d+1)}}\\
&+\Delta^{\frac{1}{\mu_d}}L^{1-\frac{2}{\mu_d}}B^{\frac{1}{\mu_d}}\Bigr].
\end{aligned}
\end{equation}
\end{lemma}

If $x_{j,d+1}=0$ for all $j$, we may take $D'=0$. Then $\Delta=1$, and the inner sums in \eqref{eq:ls} are the values of $W_b$ at the points $(x_{j,1},\dots,x_{j,d})$. If, in addition, the points are so well separated that $B=L$, then the right-hand side of \eqref{eq:ls} is $J^{1+\varepsilon}(D^{1/(d(d+1))}L^{1-1/(d(d+1))}+L^{1-1/\mu_d})$.

\begin{proof}
For complex $b=(b_r)$ and $x\in\mathbb{R}^{d+1}$ let $F_b(x)=\sum_{r=1}^{J}b_r\,e(\sum_{i=1}^{d+1}x_ir^i)$. This function has period $1$ in each of the first $d-1$ variables, and we regard it as a function on $G=\mathbb{T}^{d-1}\times\mathbb{R}^{2}$. Let
\[
\varphi(z)=\prod_{i=1}^{d+1}2J^i\,\mathbf{1}_{[-1/(4J^i),\,1/(4J^i)]}(z_i)\qquad(z\in G),
\]
and let $\Phi(x)=\sum_{j}a_j\varphi(x-x_j)$. For $1\le r\le J$,
\[
\int_{G}\varphi(z)\,e\Bigl(\sum_{i=1}^{d+1}z_ir^i\Bigr)\,dz=\prod_{i=1}^{d+1}\frac{\sin(\pi r^i/(2J^i))}{\pi r^i/(2J^i)}\in\bigl[(2/\pi)^{d+1},1\bigr].
\]
Dividing $b_r$ by this number we obtain coefficients $b'_r$ with $|b'_r|\le(\pi/2)^{d+1}$ such that
\[
\sum_{j=1}^{L}a_jF_b(x_j)=\int_{G}\Phi(x)F_{b'}(x)\,dx .
\]
We have $\|\Phi\|_1\le L$. The integral of $\varphi(x-x_j)\varphi(x-x_{j'})$ vanishes unless the pair $(j,j')$ is counted in $B$, and it never exceeds $\|\varphi\|_2^2=\prod_i2J^i$; hence $\|\Phi\|_2^2\le 2^{d+1}BJ^{(d+1)(d+2)/2}$. Moreover, $\Phi$ is supported in the set $\mathcal{Y}$ of all $x\in G$ with $|x_d|\le D+\frac{1}{4J^d}$ and $|x_{d+1}|\le D'+\frac{1}{4J^{d+1}}$. Writing $|F_{b'}(x)|$ as the integral of $\mathbf{1}_{\{|F_{b'}(x)|>\lambda\}}$ over $\lambda>0$ and applying the Cauchy--Schwarz inequality for each $\lambda$, we get
\[
\Bigl|\int_{G}\Phi F_{b'}\Bigr|\le\int_0^\infty\min\Bigl\{L,\ \|\Phi\|_2\,\bigl|\{x\in\mathcal{Y}:\ |F_{b'}(x)|>\lambda\}\bigr|^{1/2}\Bigr\}\,d\lambda .
\]
For fixed $x_{d+1}$ the factors $e(x_{d+1}r^{d+1})$ are absorbed into the coefficients, and $F_{b'}$ becomes a sum of the form \eqref{eq:W} in the remaining variables, whose coefficients have modulus at most $(\pi/2)^{d+1}$. The variable $x_d$ runs over an interval whose length $2D+\frac{1}{2J^d}$ lies between $J^{-(d-1)}$ and $3D\le 1$, and the variable $x_{d+1}$ runs over an interval of length $2D'+\frac{1}{2J^{d+1}}\le 2\Delta J^{-(d+1)}$. Hence Proposition~\ref{prop:levelset}, followed by an integration over $x_{d+1}$, gives
\[
\bigl|\{x\in\mathcal{Y}:\ |F_{b'}(x)|>\lambda\}\bigr|\lesssim_{d,\varepsilon}J^{\varepsilon-(d+1)(d+2)/2}\,\Delta\,\Bigl(D\Bigl(\frac J\lambda\Bigr)^{d(d+1)}+\Bigl(\frac J\lambda\Bigr)^{\mu_d}\Bigr).
\]
Therefore the integrand above is
\[
\lesssim_{d,\varepsilon}\min\Bigl\{L,\ (B\Delta DJ^{\varepsilon})^{1/2}\Bigl(\frac J\lambda\Bigr)^{d(d+1)/2}\Bigr\}+\min\Bigl\{L,\ (B\Delta J^{\varepsilon})^{1/2}\Bigl(\frac J\lambda\Bigr)^{\mu_d/2}\Bigr\}.
\]
Finally, $\int_0^\infty\min\{L,A(J/\lambda)^{m}\}\,d\lambda=\frac{m}{m-1}\,JL^{1-1/m}A^{1/m}$ for $A>0$ and $m>1$. With $m=d(d+1)/2$ and $m=\mu_d/2$ this gives \eqref{eq:ls}.
\end{proof}

\subsection{Short shifts}\label{subsec:shifts}

The next lemma expresses exponential sums by bilinear forms of the type considered in Lemma~\ref{lem:ls}. The Weyl sums arise from Taylor expansions of the phases around short shifts of the variable of summation, as in Vinogradov's method; see \cite{HB} and \cite[Section~3.2]{CGZZ}. The lemma is stated for a family of phases, since this is what Section~\ref{sec:lattice} requires.

\begin{lemma}\label{lem:shift}
Let $d\ge 1$ and $J\ge 1$ be integers, and let $C_1>0$. Let $\mathcal{P}$ be a finite set. For each $h\in\mathcal{P}$ let $I_h\subset\mathbb{R}$ be a bounded interval, and let $f_h:I_h\to\mathbb{R}$ be a function with $d+1$ continuous derivatives such that $|f_h^{(d+1)}|\le C_1J^{-(d+1)}$ on $I_h$. Let $\mathcal{N}$ be the set of all pairs $(h,n)$ with $h\in\mathcal{P}$ and $n\in\mathbb{Z}$ such that $n\in I_h$ and $n+J\in I_h$, and put $A_i(h,n)=f_h^{(i)}(n)/i!$ for $(h,n)\in\mathcal{N}$ and $1\le i\le d$. Then, for all complex numbers $c_h$ with $|c_h|\le 1$,
\begin{equation}\label{eq:shift}
\begin{aligned}
&\Bigl|\sum_{h\in\mathcal{P}}c_h\sum_{n\in I_h\cap\mathbb{Z}}e\bigl(f_h(n)\bigr)\Bigr|\\
&\qquad\le J\,\#\mathcal{P}+\frac{1+2\pi C_1}{J}\,\sup\,\Bigl|\sum_{(h,n)\in\mathcal{N}}a_{h,n}\sum_{r=1}^{j}e\Bigl(\sum_{i=1}^{d}A_i(h,n)\,r^i\Bigr)\Bigr|,
\end{aligned}
\end{equation}
where the supremum is over all integers $1\le j\le J$ and all complex numbers $a_{h,n}$ with $|a_{h,n}|\le 1$.
\end{lemma}

\begin{proof}
Fix $h$, and let $\mathcal{N}_h$ be the set of all $n$ with $(h,n)\in\mathcal{N}$. If $1\le r\le J$, then $\mathcal{N}_h+r$ is contained in $I_h\cap\mathbb{Z}$, and it omits at most $J$ elements of this set. Hence
\[
\Bigl|\sum_{n\in I_h\cap\mathbb{Z}}e\bigl(f_h(n)\bigr)-\frac1J\sum_{n\in\mathcal{N}_h}\sum_{r=1}^{J}e\bigl(f_h(n+r)\bigr)\Bigr|\le J .
\]
Now let $(h,n)\in\mathcal{N}$. For real $x$ with $0\le x\le J$ the point $n+x$ lies in $I_h$, and we write
\[
f_h(n+x)=f_h(n)+\sum_{i=1}^{d}A_i(h,n)\,x^i+R_{h,n}(x).
\]
Taylor's formula for $f_h'$ gives $|R_{h,n}'(x)|\le C_1J^{-(d+1)}x^d/d!\le C_1J^{-1}$. Let $w_j(h,n)$ denote the inner sum in \eqref{eq:shift}. By partial summation,
\[
\begin{aligned}
\sum_{r=1}^{J}e\bigl(f_h(n+r)\bigr)=e\bigl(f_h(n)\bigr)\Bigl[&e\bigl(R_{h,n}(J)\bigr)\,w_J(h,n)\\
&-\sum_{j=1}^{J-1}\bigl(e(R_{h,n}(j+1))-e(R_{h,n}(j))\bigr)\,w_j(h,n)\Bigr],
\end{aligned}
\]
and $|e(R_{h,n}(j+1))-e(R_{h,n}(j))|\le 2\pi C_1J^{-1}$. We multiply by $c_h$ and sum over $(h,n)\in\mathcal{N}$. For each $j$ the result is a sum $\sum a_{h,n}w_j(h,n)$ over $\mathcal{N}$, whose coefficients have modulus at most $1$ if $j=J$ and at most $2\pi C_1J^{-1}$ if $j<J$. This gives \eqref{eq:shift}.
\end{proof}

\section{Exponential sums in one variable}\label{sec:onevar}

In this section we prove Theorem~\ref{thm:onevar} and Corollary~\ref{cor:pairs}.

\subsection{Proof of Theorem~\ref{thm:onevar}}

Let $d$, $\mathcal{K}$, $c$, $C$, $N$, $t$, $I$ and $f$ be as in Theorem~\ref{thm:onevar}. In this subsection all implied constants may depend on $d$, $\mathcal{K}$, $c$ and $C$. The estimate \eqref{eq:onevar} is trivial if $N$ is bounded, so we may assume that $N$ is sufficiently large. Since $\mathcal{K}$ is compact, there is a number $\eta>0$, which depends only on $d$ and $\mathcal{K}$, such that $\frac{d+1}{2}+\eta\le t\le d-\frac1d-\eta$. We put
\[
J=\bigl\lfloor N^{1-\frac{t}{d+1}}\bigr\rfloor,\qquad D=N^{t-d}.
\]
Since $\frac1d<1-\frac{t}{d+1}<\frac12$, we have $\frac12N^{1/d}\le J\le N^{1/2}$, and $J^{d+1}\asymp N^{d+1-t}$. Moreover,
\begin{equation}\label{eq:DJ}
DJ\asymp N^{\frac{d}{d+1}(t-d+\frac1d)}\le N^{-\eta/2},\qquad DJ^{d-1}\asymp N^{\frac{2}{d+1}(t-\frac{d+1}{2})}\ge N^{\frac{2\eta}{d+1}}.
\end{equation}
These two relations are the reason for the restriction of $t$ to the interval $(\frac{d+1}{2},\,d-\frac1d)$. For $n\in I\cap\mathbb{Z}$ and $1\le i\le d$ let $A_i(n)=f^{(i)}(n)/i!$. The next lemma shows that the two Taylor coefficients of the highest degrees separate the integers in $I$.

\begin{lemma}\label{lem:sep}
Let $N$ be sufficiently large. If $n,n'\in I\cap\mathbb{Z}$ satisfy
\[
|A_d(n)-A_d(n')|\le J^{-d}\qquad\text{and}\qquad\|A_{d-1}(n)-A_{d-1}(n')\|\le J^{-(d-1)},
\]
then $n=n'$.
\end{lemma}

\begin{proof}
The derivative of $A_d$ is $f^{(d+1)}/d!$. It is continuous and has no zeros in $I$, so that it does not change sign, and its modulus is $\gtrsim N^{t-d-1}$. By the mean value theorem, the first condition gives
\[
|n-n'|\lesssim N^{d+1-t}J^{-d}\lesssim J .
\]
In the same way, the derivative of $A_{d-1}$ is $f^{(d)}/(d-1)!$, which does not change sign and whose modulus is $\asymp D$. Hence $|A_{d-1}(n)-A_{d-1}(n')|\asymp D\,|n-n'|\lesssim DJ$. By \eqref{eq:DJ} this is less than $\frac12$ if $N$ is large. Then the distance from $A_{d-1}(n)-A_{d-1}(n')$ to the nearest integer is equal to its modulus, and the second condition gives $|n-n'|\lesssim J^{-(d-1)}/D$. By \eqref{eq:DJ} the right-hand side is less than $1$ if $N$ is large, and then $n=n'$.
\end{proof}

\begin{proof}[Proof of Theorem~\ref{thm:onevar}]
By the choice of $J$ we have $|f^{(d+1)}|\le CN^{t-d-1}\le CJ^{-(d+1)}$ on $I$. Hence Lemma~\ref{lem:shift}, applied to the single function $f$ with $C_1=C$, gives
\[
\Bigl|\sum_{n\in I\cap\mathbb{Z}}e\bigl(f(n)\bigr)\Bigr|\le J+\frac{1+2\pi C}{J}\,\sup\,\Bigl|\sum_{n\in\mathcal{N}}a_n\sum_{r=1}^{j}e\Bigl(\sum_{i=1}^{d}A_i(n)\,r^i\Bigr)\Bigr|,
\]
where $\mathcal{N}$ is the set of all integers $n$ with $n\in I$ and $n+J\in I$, and the supremum is over all integers $1\le j\le J$ and all complex numbers $a_n$ with $|a_n|\le 1$. The inner sum is the value of $W_b$ at the point $(A_1(n),\dots,A_d(n))$, for $b_r=\mathbf{1}_{\{r\le j\}}$.

We apply Lemma~\ref{lem:ls} to the points $(A_1(n),\dots,A_d(n),0)$ with $n\in\mathcal{N}$, and with $D'=0$. Their number is $L=\#\mathcal{N}\le N+1$. We have $|A_d(n)|\le CD$, and $J^{-(d-1)}\le CD\le\frac14$ if $N$ is large; this follows from \eqref{eq:DJ} and from $D\le N^{-1/d}$. Thus the lemma can be applied with $CD$ in place of $D$, and with $\Delta=1$. By Lemma~\ref{lem:sep}, the only pairs $(n,n')$ that are counted in Lemma~\ref{lem:ls} are those with $n=n'$, so that $B=L$. Therefore \eqref{eq:ls} bounds the supremum by
\[
C_\varepsilon J^{1+\varepsilon}\Bigl(D^{\frac{1}{d(d+1)}}N^{1-\frac{1}{d(d+1)}}+N^{1-\frac{1}{\mu_d}}\Bigr).
\]
Since $D=N^{t-d}$, the first term in the parentheses is equal to $N^{1-\frac{d+1-t}{d(d+1)}}$. Moreover, $J\le N^{1/2}\le N^{1-\frac{1}{\mu_d}}$. This proves \eqref{eq:onevar}.
\end{proof}

\subsection{Exponent pairs}\label{subsec:pairs}

An exponent pair is a bound for the sums $\sum e(f(n))$ that holds for all sizes of the phase $f$. We recall the definition from \cite[Definition~5.1]{DB}, in the equivalent form of \cite[Lemma~5.2]{DB}.

\begin{definition}\label{def:pairs}
Let $k$ and $\ell$ be real numbers with $0\le k\le\frac12\le\ell\le 1$ and $k+\ell\le 1$. The pair $(k,\ell)$ is called an exponent pair if for all $\varepsilon>0$ and $\sigma>0$ there are $\delta>0$, an integer $P\ge 1$ and $C\ge 1$ such that
\begin{equation}\label{eq:pairbound}
\Bigl|\sum_{n\in I\cap\mathbb{Z}}e\Bigl(N^{t}F\Bigl(\frac nN\Bigr)\Bigr)\Bigr|\le C\,N^{\ell+(t-1)k+\varepsilon t}
\end{equation}
for all real numbers $N\ge 2$ and $t\ge 1$, all intervals $I\subseteq[N,2N]$, and all smooth functions $F:[1,2]\to\mathbb{R}$ with
\begin{equation}\label{eq:model}
\Bigl|\frac{\mathrm{d}^i}{\mathrm{d}u^i}\bigl(F'(u)-u^{-\sigma}\bigr)\Bigr|\le\delta\qquad(0\le i\le P,\ u\in[1,2]).
\end{equation}
\end{definition}

In \cite{DB} the phase is written as $TF(n/N)$ with $T\ge N\ge 1$, and the bound as $C(T/N)^{k+\varepsilon}N^{\ell+\varepsilon}$. This is \eqref{eq:pairbound} with $T=N^{t}$; for $N<2$ the sum has at most two terms. The condition \eqref{eq:model} says that $F'$ is close to the function $u^{-\sigma}$, together with its first $P$ derivatives. The bound \eqref{eq:pairbound} is trivial if $\ell+(t-1)k\ge 1$, because the sum has at most $N+1$ terms.

Let $\sigma>0$, and let $P\ge 1$ be an integer. For $i\ge 0$ the derivative of order $i$ of $u^{-\sigma}$ has the modulus $\sigma(\sigma+1)\cdots(\sigma+i-1)\,u^{-\sigma-i}$, which lies between two positive constants for $u\in[1,2]$. Hence there are numbers $\delta_0>0$ and $0<c'\le C'$, which depend only on $\sigma$ and $P$, with the following property: if $F$ satisfies \eqref{eq:model} with some $\delta\le\delta_0$, and if $N\ge 2$ and $t\ge 1$, then the function $f(x)=N^{t}F(x/N)$ satisfies
\begin{equation}\label{eq:fsize}
c'N^{t-i}\le|f^{(i)}(x)|\le C'N^{t-i}\qquad(x\in[N,2N],\ 1\le i\le P+1).
\end{equation}
Thus Theorem~\ref{thm:onevar} can be applied to $f$ if $P\ge d$. We combine it with the estimate of Heath-Brown~\cite{HB}, which we use in the following form.

\begin{lemma}\label{lem:hb}
Let $m\ge 3$ be an integer, and let $0<c\le C$. Let $N\ge 2$ and $t\ge 1$, let $I\subseteq[N,2N]$ be an interval, and let $f:I\to\mathbb{R}$ be a function with $m$ continuous derivatives such that $cN^{t-m}\le|f^{(m)}(x)|\le CN^{t-m}$ for $x\in I$. Then
\[
\sum_{n\in I\cap\mathbb{Z}}e\bigl(f(n)\bigr)\lesssim_{m,c,C,\varepsilon}N^{1+\varepsilon}\Bigl(N^{-\frac{m-t}{m(m-1)}}+N^{-\frac{1}{m(m-1)}}+N^{-\frac{2t}{m^2(m-1)}}\Bigr).
\]
In particular, the sum is $O_{m,c,C,\varepsilon}(N^{1-\frac{1}{m(m-1)}+\varepsilon})$ if $\frac m2\le t\le m-1$.
\end{lemma}

\begin{proof}
This is \cite[Theorem~1]{HB} with $\lambda_m=cN^{t-m}$ and $A=C/c$, applied to $f$ or to $-f$ after a translation of the variable. The three terms of the bound in \cite{HB} are $\lambda_m^{1/(m(m-1))}$, $N^{-1/(m(m-1))}$ and $N^{-2/(m(m-1))}\lambda_m^{-2/(m^2(m-1))}$. In \cite{HB} the sum is over the integers $1\le n\le N$. Since the bound of that theorem increases with $N$, it can also be used for the sum over $I\cap\mathbb{Z}$, which has at most $N+1$ terms. If $\frac m2\le t\le m-1$, the second term is the largest.
\end{proof}

For integers $j\ge 4$ let
\[
t_j=j-\frac1j,\qquad w_j=1-\frac{1}{j(j+1)}.
\]
By Lemma~\ref{lem:hb} with $m=j+1$, the sums under consideration are $O(N^{w_j+\varepsilon})$ for $\frac{j+1}{2}\le t\le j$. Theorem~\ref{thm:onevar} with $d=j$ gives a smaller exponent when $t$ is slightly smaller than $t_j$, but it does not apply for $t\ge t_j$. In this sense the points $(t_j,w_j)$ are the corners of the bounds that are at our disposal, and an exponent pair corresponds to a line that passes above these corners. For $d\ge 4$ let
\[
\Lambda_d(t)=\ell_d+(t-1)k_d,
\]
with $k_d$ and $\ell_d$ as in \eqref{eq:kl}. Since $w_{d+1}-w_d=\frac{2}{d(d+1)(d+2)}$ and $t_{d+1}-t_d=\frac{d^2+d+1}{d(d+1)}$, a short computation shows that $\Lambda_d$ is the linear function with
\begin{equation}\label{eq:chord}
\Lambda_d(t_d)=w_d,\qquad\Lambda_d(t_{d+1})=w_{d+1}.
\end{equation}
The next lemma is the main step. It shows that a line through two consecutive corners, or a line with similar properties, is an upper bound for the exponents between them.

\begin{lemma}\label{lem:chord}
Let $j\ge 3$ be an integer, and let $\Lambda$ be a linear function whose slope lies between $0$ and $\frac{1}{(j+1)(j+2)}$. Suppose that $\Lambda(t_{j+1})\ge w_{j+1}$, and that $\Lambda(t')\ge 1-\frac{1}{\mu_{j+1}}$ for some $t'$ with $\frac{j+1}{2}\le t'\le t_{j+1}$. Let $0<\eta<\frac12$ and $0<c\le C$. Let $N\ge 2$ and $t'\le t\le t_{j+1}$, let $I\subseteq[N,2N]$ be an interval, and let $f:I\to\mathbb{R}$ be a function with $j+2$ continuous derivatives such that
\[
cN^{t-i}\le|f^{(i)}(x)|\le CN^{t-i}\qquad(x\in I,\ i\in\{j+1,j+2\}).
\]
Then
\[
\sum_{n\in I\cap\mathbb{Z}}e\bigl(f(n)\bigr)\lesssim_{j,\eta,c,C}N^{\Lambda(t)+\eta}.
\]
\end{lemma}

Note that $1-\frac{1}{\mu_{j+1}}=w_j$ if $j\ge 4$, and that $1-\frac{1}{\mu_{j+1}}=\frac{13}{14}$ if $j=3$.

\begin{proof}
Since the slope of $\Lambda$ is not negative, $\Lambda(t)\ge\Lambda(t')\ge 1-\frac{1}{\mu_{j+1}}$. We distinguish three cases, and we use Lemma~\ref{lem:hb} and Theorem~\ref{thm:onevar} with $\varepsilon=\eta/2$.

If $t\le j$, Lemma~\ref{lem:hb} with $m=j+1$ bounds the sum by $O(N^{1-\frac{1}{j(j+1)}+\eta/2})$, because $\frac{j+1}{2}\le t\le j$; note that $j(j+1)\le\mu_{j+1}$.

If $j\le t\le t_{j+1}-\frac\eta2$, we apply Theorem~\ref{thm:onevar} with $d=j+1$ and $\mathcal{K}=[j,\,t_{j+1}-\frac\eta2]$, which is a compact subinterval of $(\frac{j+2}{2},\,t_{j+1})$. It bounds the sum by $O(N^{1-\frac{1}{\mu_{j+1}}+\eta/2}+N^{\Gamma(t)+\eta/2})$, where
\[
\Gamma(t)=1-\frac{j+2-t}{(j+1)(j+2)}.
\]
The slope of the linear function $\Gamma$ is $\frac{1}{(j+1)(j+2)}$, which is at least the slope of $\Lambda$, and
\[
\Gamma(t_{j+1})=1-\frac{1}{(j+1)^2}<w_{j+1}\le\Lambda(t_{j+1}).
\]
Hence $\Gamma(t)\le\Lambda(t)$ for all $t\le t_{j+1}$.

If $t_{j+1}-\frac\eta2\le t\le t_{j+1}$, Lemma~\ref{lem:hb} with $m=j+2$ bounds the sum by $O(N^{w_{j+1}+\eta/2})$, because $\frac{j+2}{2}\le t\le j+1$. Since the slope of $\Lambda$ is less than $1$, we have $w_{j+1}\le\Lambda(t_{j+1})\le\Lambda(t)+\frac\eta2$.
\end{proof}

\begin{lemma}\label{lem:ranges}
Let $j_0\ge 4$ be an integer, and let $0<c\le C$. Let $N\ge 2$ and $1\le t\le t_{j_0+1}$, let $I\subseteq[N,2N]$ be an interval, and let $f:I\to\mathbb{R}$ be a function with $j_0+2$ continuous derivatives such that $cN^{t-i}\le|f^{(i)}(x)|\le CN^{t-i}$ for $x\in I$ and $3\le i\le j_0+2$. Then the sum $\sum_{n\in I\cap\mathbb{Z}}e(f(n))$ is
\begin{enumerate}
\item[(i)] $O(N^{8/9+\varepsilon})$ if $1\le t\le 2$, and $O(N^{11/12+\varepsilon})$ if $2\le t\le 3$;
\item[(ii)] $O(N^{\Lambda_4(t)+\varepsilon})$ if $\frac{12}{5}\le t\le t_4$;
\item[(iii)] $O(N^{\Lambda_j(t)+\varepsilon})$ if $t_j\le t\le t_{j+1}$ for an integer $4\le j\le j_0$.
\end{enumerate}
The implied constants depend only on $j_0$, $c$, $C$ and $\varepsilon$.
\end{lemma}

\begin{proof}
(i) For $1\le t\le 2$ the three terms in Lemma~\ref{lem:hb} with $m=3$ are $N^{-(3-t)/6}$, $N^{-1/6}$ and $N^{-t/9}$, so that the sum is $O(N^{8/9+\varepsilon})$. For $2\le t\le 3$ the lemma with $m=4$ gives $O(N^{11/12+\varepsilon})$.

(ii) We may assume that $\varepsilon<\frac12$. We apply Lemma~\ref{lem:chord} with $j=3$, $\Lambda=\Lambda_4$, $t'=\frac{12}{5}$ and $\eta=\varepsilon$. Its hypotheses hold: the slope of $\Lambda_4$ is $k_4=\frac{1}{63}<\frac{1}{20}$, we have $\Lambda_4(t_4)=w_4$ by \eqref{eq:chord}, and $\Lambda_4(\frac{12}{5})=\ell_4+\frac75k_4=\frac{13}{14}=1-\frac{1}{\mu_4}$.

(iii) We apply Lemma~\ref{lem:chord} with $\Lambda=\Lambda_j$, $t'=t_j$ and $\eta=\varepsilon$. Its hypotheses hold by \eqref{eq:chord}, and because $k_j\le\frac{1}{(j+1)(j+2)}$; the latter inequality is equivalent to $2(j+1)\le j^2+j+1$.
\end{proof}

\begin{proof}[Proof of Corollary~\ref{cor:pairs}]
Let $d\ge 4$. It is clear from \eqref{eq:kl} that $0<k_d<\frac12<\ell_d<1$ and $k_d+\ell_d<1$. By \eqref{eq:hbpairs}, $p_{d+1}=\frac{2}{d^2(d+3)}$ and $q_{d+1}=1-\frac{3d+1}{d(d+1)(d+3)}$. The inequalities $k_d<p_{d+1}$ and $\ell_d<q_{d+1}$ therefore follow from
\[
\begin{aligned}
(d+2)(d^2+d+1)-d^2(d+3)&=3d+2>0,\\
(3d^2+3d-1)\,d(d+3)-(3d+1)(d+2)(d^2+d+1)&=2(d^3-2d^2-6d-1)>0.
\end{aligned}
\]

It remains to show that $(k_d,\ell_d)$ is an exponent pair. Let $\varepsilon>0$ and $\sigma>0$. We choose an integer $j_0\ge d$ with $\Lambda_d(t_{j_0+1})\ge 1$, and we take $P=j_0+1$ and $\delta=\delta_0$, where $\delta_0$ is the number that belongs to $\sigma$ and $P$ in \eqref{eq:fsize}. Let $N$, $t$, $I$ and $F$ be as in Definition~\ref{def:pairs}, and let $f(x)=N^{t}F(x/N)$. If $t\ge t_{j_0+1}$, then $\Lambda_d(t)\ge 1$, and \eqref{eq:pairbound} is trivial. Let $1\le t\le t_{j_0+1}$. By \eqref{eq:fsize}, Lemma~\ref{lem:ranges} applies to $f$, with constants that depend only on $\sigma$, $j_0$ and $\varepsilon$. Since $N^{\varepsilon}\le N^{\varepsilon t}$, it suffices to show that for every $t$ in this range one of the exponents in that lemma is at most $\Lambda_d(t)$.

For this we compare the functions $\Lambda_j$ with each other. Let $j\ge 4$. By \eqref{eq:chord}, both $\Lambda_j$ and $\Lambda_{j+1}$ take the value $w_{j+1}$ at $t_{j+1}$, so that
\[
\Lambda_{j+1}(t)-\Lambda_j(t)=(k_{j+1}-k_j)(t-t_{j+1}),
\]
and $k_{j+1}<k_j$. Hence $\Lambda_j\le\Lambda_{j+1}$ on $(-\infty,t_{j+1}]$, and $\Lambda_{j+1}\le\Lambda_j$ on $[t_{j+1},\infty)$. Since the numbers $t_j$ increase with $j$, it follows by induction that
\begin{equation}\label{eq:lines}
\Lambda_j(t)\le\Lambda_d(t)\quad\text{if}\quad 4\le j\le d,\ t\le t_{j+1},\qquad\text{or if}\quad j\ge d,\ t\ge t_j .
\end{equation}
In particular, $\Lambda_j\le\Lambda_d$ on $[t_j,t_{j+1}]$ for every $j\ge 4$, and $\Lambda_4\le\Lambda_d$ on $[1,t_4]$. This settles the range $\frac{12}{5}\le t\le t_{j_0+1}$, by the parts (ii) and (iii) of Lemma~\ref{lem:ranges}. Finally, $\frac89<\ell_4=\Lambda_4(1)$ and $\frac{11}{12}<\frac{83}{90}=\Lambda_4(2)$. Hence the exponents in part (i) are at most $\Lambda_4(t)\le\Lambda_d(t)$ for $1\le t\le 2$ and for $2\le t\le\frac{12}{5}$, respectively. This completes the proof.
\end{proof}

The same method gives one more exponent pair. It lies between the pair $(k_4,\ell_4)$ and the exponent pair
\[
(k',\ell')=\Bigl(\frac{89}{1282},\ \frac{997}{1282}\Bigr)
\]
of Tao, Trudgian and Yang~\cite{TTY}; see also \cite[Theorem~5.22]{DB}.

\begin{proposition}\label{prop:star}
The pair
\[
(k_*,\ell_*)=\Bigl(\frac{89}{2395},\ \frac{4061}{4790}\Bigr)
\]
is an exponent pair.
\end{proposition}

In the list of \cite[Chapter~5]{DB}, the pair $(k',\ell')$ is a vertex of the convex hull of the known exponent pairs. As we shall see in Section~\ref{subsec:hull}, $(k_*,\ell_*)$ is a vertex of the convex hull of the known exponent pairs and the pairs of Corollary~\ref{cor:pairs}. In particular, it is not a convex combination of $(k',\ell')$ and the pairs $(k_d,\ell_d)$.

\begin{proof}
Let $\Lambda_*(t)=\ell_*+(t-1)k_*$ and $t_*=\frac{1977}{623}$. A direct computation gives
\[
\Lambda_*(t_*)=\frac{13}{14}=1-\frac{1}{\mu_4},\qquad\Lambda_*(t_4)=\frac{19}{20}=w_4,\qquad\ell'+(t_*-1)k'=\frac{13}{14}.
\]
Thus $\Lambda_*$ is the line through the corner $(t_4,w_4)$ and the point at which the bound of the pair $(k',\ell')$ reaches the exponent $\frac{13}{14}$.

Clearly $0<k_*<\frac12<\ell_*<1$ and $k_*+\ell_*<1$. Let $\varepsilon>0$ and $\sigma>0$, where we may assume that $\varepsilon<\frac12$. We choose an integer $j_0\ge 4$ with $\Lambda_4(t_{j_0+1})\ge 1$. Since $(k',\ell')$ is an exponent pair, Definition~\ref{def:pairs} provides numbers $\delta_1$ and $P_1$ for this pair and for $\varepsilon$ and $\sigma$. We take $P=\max(P_1,j_0+1)$ and $\delta=\min(\delta_1,\delta_0)$, with $\delta_0$ as in \eqref{eq:fsize}. Let $N$, $t$, $I$ and $F$ be as in Definition~\ref{def:pairs}, and let $f(x)=N^{t}F(x/N)$.

If $1\le t\le t_*$, the sum in \eqref{eq:pairbound} is $O(N^{\ell'+(t-1)k'+\varepsilon t})$. The linear function $\ell'+(t-1)k'$ agrees with $\Lambda_*$ at $t_*$, and its slope $k'$ is larger than $k_*$. Hence $\ell'+(t-1)k'\le\Lambda_*(t)$ for $t\le t_*$.

If $t_*\le t\le t_4$, we apply Lemma~\ref{lem:chord} with $j=3$, $\Lambda=\Lambda_*$, $t'=t_*$ and $\eta=\varepsilon$; note that $k_*<\frac{1}{20}$ and $2\le t_*\le t_4$. It bounds the sum by $O(N^{\Lambda_*(t)+\varepsilon})$.

Let $t\ge t_4$. Since $\Lambda_*(t_4)=\Lambda_4(t_4)$ and $k_*>k_4$, we have $\Lambda_4(t)\le\Lambda_*(t)$. If $t\ge t_{j_0+1}$, then $\Lambda_*(t)\ge 1$, and \eqref{eq:pairbound} is trivial. Otherwise $t_j\le t\le t_{j+1}$ for some $4\le j\le j_0$, and by Lemma~\ref{lem:ranges} and \eqref{eq:lines} with $d=4$ the sum is $O(N^{\Lambda_4(t)+\varepsilon})$.

In all cases the sum is $O(N^{\Lambda_*(t)+\varepsilon t})$, with a constant that depends only on $\varepsilon$ and $\sigma$. Hence $(k_*,\ell_*)$ is an exponent pair.
\end{proof}

\subsection{Comparison with the known exponent pairs}\label{subsec:hull}

\begin{table}[t]
\caption{The new exponent pairs that are closest to the vertices of $\mathcal{H}$. The fourth column gives the boundary of $\mathcal{H}$ at the first coordinate of the pair, and the last column the two vertices of $\mathcal{H}$ between which this boundary point lies. All numbers are rounded to five decimal places.}\label{tab:hull}
\small
\renewcommand{\arraystretch}{1.2}
\begin{tabular}{@{}llllll@{}}
\hline
pair & $k$ & $\ell$ & $\ell_{\mathcal{H}}(k)$ & $\ell_{\mathcal{H}}(k)-\ell$ & edge of $\mathcal{H}$\\
\hline
$(k_*,\ell_*)$ & $0.03716$ & $0.84781$ & $0.85203$ & $0.00422$ & $V_8$, $V_7$\\
$(k_4,\ell_4)$ & $0.01587$ & $0.90635$ & $0.91376$ & $0.00741$ & $(p_6,q_6)$, $V_9$\\
$(k_5,\ell_5)$ & $0.00922$ & $0.93164$ & $0.93669$ & $0.00505$ & $(p_7,q_7)$, $(p_6,q_6)$\\
$(k_6,\ell_6)$ & $0.00581$ & $0.94809$ & $0.95162$ & $0.00353$ & $(p_8,q_8)$, $(p_7,q_7)$\\
$(k_7,\ell_7)$ & $0.00390$ & $0.95931$ & $0.96185$ & $0.00254$ & $(p_9,q_9)$, $(p_8,q_8)$\\
$(k_8,\ell_8)$ & $0.00274$ & $0.96728$ & $0.96916$ & $0.00189$ & $(p_{10},q_{10})$, $(p_9,q_9)$\\
\hline
\end{tabular}
\end{table}

The exponent pairs that are known so far are described in \cite[Chapter~5]{DB} by a table of the vertices of their convex hull. We write $V_0,\dots,V_9$ for the pairs in the rows $0,\dots,9$ of that table. Thus $V_0=(\frac{13}{84},\frac{55}{84})$ is the exponent pair of Bourgain, and $V_4=(k',\ell')$. The other rows contain the pairs $(p_m,q_m)$ with $m\ge 6$, and the images of all these pairs under the map
\[
B(k,\ell)=\Bigl(\ell-\frac12,\ k+\frac12\Bigr),
\]
which describes the van der Corput $B$-process \cite[Proposition~5.9]{DB}. Let $\mathcal{H}$ be the convex hull of the pairs in this table and of the exponent pairs $(0,1)$ and $(\frac12,\frac12)$; according to \cite[Chapter~5]{DB}, $\mathcal{H}$ is the set of the exponent pairs that are known so far. Let $\mathcal{H}'$ be the convex hull of $\mathcal{H}$, of the pairs $(k_d,\ell_d)$ with $d\ge 4$, of the pair $(k_*,\ell_*)$, and of the images of these pairs under $B$. Since the set of exponent pairs is convex \cite[Corollary~5.4]{DB} and closed under $B$, every element of $\mathcal{H}'$ is an exponent pair.

A direct computation shows that the vertices of $\mathcal{H}'$ are $(0,1)$, $(\frac12,\frac12)$, the pairs $(k_d,\ell_d)$ with $d\ge 4$, the pair $(k_*,\ell_*)$, the pairs $V_0,\dots,V_5$, and the images of all these pairs under $B$. In particular, the pairs $(k_d,\ell_d)$ and $(k_*,\ell_*)$ do not belong to $\mathcal{H}$, whereas $V_6,\dots,V_9$ and the pairs $(p_m,q_m)$ with $m\ge 6$ are no longer vertices.

Table~\ref{tab:hull} shows how far the new pairs that are closest to the old vertices lie outside $\mathcal{H}$. Here, for $0\le k\le\frac{13}{84}$, $\ell_{\mathcal{H}}(k)$ is the smallest number $\ell$ with $(k,\ell)\in\mathcal{H}$, so that a pair $(k,\ell)$ with $k\le\frac{13}{84}$ and $k+\ell\le 1$ lies outside $\mathcal{H}$ if and only if $\ell<\ell_{\mathcal{H}}(k)$.

\section{The circle and divisor problems}\label{sec:lattice}

In this section we prove Theorem~\ref{thm:main}. We use the notation of \cite{CGZZ}, which is introduced in Section~\ref{subsec:reduction}. That subsection also contains the reduction of the two remainder terms to exponential sums and the statement of the main estimate, Theorem~\ref{thm:rec}. The proof of Theorem~\ref{thm:rec} occupies Sections~\ref{subsec:classical}--\ref{subsec:large}. Classical bounds cover a part of the range of the parameters (Section~\ref{subsec:classical}). The rest is divided into a middle range, which is treated in Sections~\ref{subsec:middle}--\ref{subsec:close} with Lemmas~\ref{lem:shift} and~\ref{lem:ls}, and a second range, which is treated in Section~\ref{subsec:large} with the method of \cite[Section~5]{CGZZ}. Theorem~\ref{thm:main} is deduced in Section~\ref{subsec:proofmain}.

\emph{The use of \cite{CGZZ}.} The proof of Theorem~\ref{thm:rec} follows the proof of \cite[Theorem~1.2]{CGZZ}, and several steps are taken from \cite{CGZZ}. They are of three kinds.
\begin{enumerate}
\item[(a)] Two results of \cite{CGZZ} hold for all sums \eqref{eq:S}, and we use them as they are stated. These are the bounds of \cite[Proposition~2.2]{CGZZ}, which follow from classical exponent pairs, and the van der Corput transformation \cite[(3.8)]{CGZZ}.
\item[(b)] Two results are proved in \cite{CGZZ} for a smaller range of the parameters, namely the count of close pairs \cite[Proposition~3.1]{CGZZ} and the estimate \cite[Lemma~5.1]{CGZZ}, which is based on the work of Watt~\cite{Watt}. We use their proofs in our range. In these proofs the parameters enter only through identities between them and through a finite list of inequalities. In Sections~\ref{subsec:close} and~\ref{subsec:large} we list these inequalities, together with the bounds that are proved in \cite{CGZZ} (Tables~\ref{tab:close} and~\ref{tab:watt}), and we prove them in our range. The proofs in \cite{CGZZ} then apply without change.
\item[(c)] The estimate \cite[Lemma~3.1]{CGZZ} for the bilinear sums is not used. It is replaced by Lemma~\ref{lem:csieve}, which is proved in Section~\ref{subsec:sieveproof} with Lemmas~\ref{lem:shift} and~\ref{lem:ls}.
\end{enumerate}
The parameters that change are listed in Table~\ref{tab:dict}. Apart from a few conditions that hold for all values of the parameters, each inequality in (b) compares two powers of $T$ whose exponents are piecewise linear functions of the logarithms of the parameters. The extreme values of these exponents are therefore attained at vertices of polygons, and the tables indicate where.

\begin{table}[ht]
\caption{The parameters of \cite{CGZZ} and the corresponding ones in this paper.}\label{tab:dict}
\small
\renewcommand{\arraystretch}{1.2}
\begin{tabular}{@{}lll@{}}
\hline
 & \cite{CGZZ} & this paper\\
\hline
exponent $\beta$ & $209/670$ & $293/940$\\
threshold $s$ for $H$ & $123/1340$ & $171/1880$\\
middle range & $T^{42/67}\le\rho\le T^{43/67}$, \cite[(2.3)]{CGZZ} & $T^{42/67}\le\rho\le T^{61/94}$, \eqref{eq:middle}\\
length $J$ of the shifts & $\lceil T^{36/67}/M\rceil$, \cite[(3.1)]{CGZZ} & $\lceil(T\rho^{11})^{1/15}/M\rceil$, \eqref{eq:J}\\
larger-$\rho$ range & $\beta\le\tau\le 24/67$ & $\beta\le\tau\le 33/94$, \eqref{eq:large}\\
bilinear sums & \cite[Lemma~3.1]{CGZZ} & Lemma~\ref{lem:csieve}, new proof\\
close pairs & \cite[Proposition~3.1]{CGZZ} & Proposition~\ref{prop:close}, proof of \cite{CGZZ}\\
Watt's method & \cite[Lemma~5.1]{CGZZ} & \eqref{eq:watt}, proof of \cite{CGZZ}\\
\hline
\end{tabular}
\end{table}

\subsection{From lattice points to exponential sums}\label{subsec:reduction}

In this subsection we recall how $E(X)$ and $\Delta(X)$ are reduced to exponential sums, and we explain the parameters that are used in the rest of this section. The reduction is classical; we follow \cite[Section~2]{CGZZ} and \cite{LY}.

\emph{Sums of the sawtooth function.} Recall that $\psi(x)=x-\lfloor x\rfloor-\frac12$. Let $\chi$ be the nontrivial character modulo~$4$. Every integer $n\ge 1$ has $4\sum_{a\mid n}\chi(a)$ representations as a sum of two squares. Hence the number of pairs $(u,v)\in\mathbb{Z}^2$ with $u^2+v^2\le X$ is $1+4\sum_{ab\le X}\chi(a)$, where $a$ and $b$ run over the positive integers. By the hyperbola method, this sum is equal to
\[
\sum_{a\le\sqrt{X}}\chi(a)\Bigl\lfloor\frac Xa\Bigr\rfloor+\sum_{b\le\sqrt{X}}\ \sum_{a\le X/b}\chi(a)-\bigl\lfloor\sqrt{X}\bigr\rfloor\sum_{a\le\sqrt{X}}\chi(a).
\]
In the first sum we insert $\lfloor x\rfloor=x-\frac12-\psi(x)$. In the second sum we use the identity
\[
\sum_{a\le x}\chi(a)=\frac12+\psi\Bigl(\frac x4+\frac14\Bigr)-\psi\Bigl(\frac x4-\frac14\Bigr)\qquad(x>0),
\]
which follows from the same formula for $\lfloor\cdot\rfloor$, because the left-hand side is equal to $\lfloor\frac{x+3}{4}\rfloor-\lfloor\frac{x+1}{4}\rfloor$. A short computation with $\sum_{a\ge 1}\chi(a)/a=\pi/4$ shows that the terms without $\psi$ contribute $\pi X+O(1)$. Therefore
\begin{equation}\label{eq:Epsi}
\begin{aligned}
E(X)={}&-4\sum_{a\le\sqrt{X}}\chi(a)\,\psi\Bigl(\frac Xa\Bigr)\\
&+4\sum_{n\le\sqrt{X}}\Bigl(\psi\Bigl(\frac{X}{4n}+\frac14\Bigr)-\psi\Bigl(\frac{X}{4n}-\frac14\Bigr)\Bigr)+O(1).
\end{aligned}
\end{equation}
For the divisor problem, the hyperbola method gives
\[
\#\{(a,b)\in\mathbb{N}^2:\ ab\le X\}=2\sum_{n\le\sqrt{X}}\Bigl\lfloor\frac Xn\Bigr\rfloor-\bigl\lfloor\sqrt{X}\bigr\rfloor^2 .
\]
The same substitution, together with the formula $\sum_{n\le y}1/n=\log y+\gamma-\psi(y)/y+O(y^{-2})$, gives
\begin{equation}\label{eq:Dpsi}
\Delta(X)=-2\sum_{n\le\sqrt{X}}\psi\Bigl(\frac Xn\Bigr)+O(1).
\end{equation}

\emph{A common form.} We split the ranges of summation in \eqref{eq:Epsi} and \eqref{eq:Dpsi} into dyadic blocks, and in the first sum of \eqref{eq:Epsi} we write $a=4n+1$ or $a=4n-1$; the term with $a=1$ is $O(1)$. Each of the resulting $O(\log X)$ sums has the form
\begin{equation}\label{eq:psiblock}
\sum_{n\in I\cap\mathbb{Z}}\psi\Bigl(\frac TM\,F\Bigl(\frac nM\Bigr)\Bigr),\qquad F(z)=\frac{1}{4(z+\alpha)}+\kappa,
\end{equation}
where $1\le M\le\sqrt{X}$, $I\subseteq[M,2M)$ is an interval, and $T$, $\alpha$ and $\kappa$ are chosen as follows. For the arguments $X/(4n\pm1)$ we take $T=X$ and $(\alpha,\kappa)=(\pm\frac{1}{4M},0)$; for the arguments $X/(4n)\pm\frac14$ we take $T=X$ and $(\alpha,\kappa)=(0,\pm\frac{M}{4X})$; and for the arguments $X/n$ in \eqref{eq:Dpsi} we take $T=4X$ and $\alpha=\kappa=0$. This works because
\[
\frac TM\,F\Bigl(\frac nM\Bigr)=\frac{T}{4(n+M\alpha)}+\frac{T\kappa}{M}.
\]
In all cases $M\le T^{1/2}$, $|\alpha|\le\frac14$ and $|\kappa|\le 1$. Thus $T$ measures the size of the problem, $M$ is the size of the summation variable, and the small shifts $\alpha$ and $\kappa$ come from the residue classes modulo $4$ in \eqref{eq:Epsi}.

\emph{Fourier expansion.} Throughout this section we put
\[
\beta=\frac{293}{940},\qquad s=\frac72\beta-1=\frac{171}{1880}.
\]
In \cite{CGZZ} these two letters denote the numbers $209/670$ and $123/1340$. The number $\beta$ is the exponent in Theorem~\ref{thm:main}, and the role of $s$ is explained at the end of this subsection. If $M<2T^{\beta}$, the sum \eqref{eq:psiblock} is trivially $O(T^{\beta})$. If $M\ge 2T^{\beta}$, we approximate $\psi$ by Vaaler's trigonometric polynomial~\cite{Vaaler} with the frequencies $1\le|h|\le MT^{-\beta}$. The error of this approximation is bounded by a nonnegative trigonometric polynomial with the same frequencies, whose constant term is $O(T^{\beta}/M)$. The coefficients of both polynomials are $O(1/|h|)$, and their total variation on each dyadic block $H\le|h|<2H$ is $O(1/H)$. We insert this into \eqref{eq:psiblock}, sum the error polynomial over $n$, and use partial summation in $h$ on each dyadic block. Since the terms with $h<0$ are the complex conjugates of those with $h>0$, this gives \cite[(2.1)]{CGZZ}
\begin{equation}\label{eq:vaaler}
\begin{aligned}
&\Bigl|\sum_{n\in I\cap\mathbb{Z}}\psi\Bigl(\frac TM\,F\Bigl(\frac nM\Bigr)\Bigr)\Bigr|\\
&\qquad\lesssim T^{\beta}+\sum_{H=2^j\le MT^{-\beta}}\frac1H\sup_{H\le H'\le 2H}\Bigl|\sum_{\substack{H\le h<H'\\ n\in I\cap\mathbb{Z}}}e\Bigl(\frac{hT}{M}\,F\Bigl(\frac nM\Bigr)\Bigr)\Bigr|.
\end{aligned}
\end{equation}
The term $T^{\beta}$ comes from the constant term of the error polynomial; this is why the frequencies are cut off at $MT^{-\beta}$.

\emph{The exponential sums.} The inner sums in \eqref{eq:vaaler} are double sums over a frequency $h\approx H$ and the original variable $n\approx M$. Following \cite[Theorem~1.2]{CGZZ}, we consider the more general sums
\begin{equation}\label{eq:S}
S=\sum_{h\approx H}g\Bigl(\frac hH\Bigr)\sum_{n\approx M}G\Bigl(\frac nM\Bigr)\,e\Bigl(\frac{hT}{M}\,F\Bigl(\frac nM\Bigr)\Bigr),
\end{equation}
with $F$ as in \eqref{eq:psiblock}, where $T\ge 2$, $H\ge 1$, $1\le M\le T^{1/2}$, $|\alpha|\le\frac14$, $|\kappa|\le 1$, and the weights $g,G:[1,2]\to\mathbb{C}$ satisfy
\[
\|g\|_{\infty}+\operatorname{Var}(g)\le 1,\qquad \|G\|_{\infty}+\operatorname{Var}(G)\le 1 .
\]
Each of the ranges $h\approx H$ and $n\approx M$ may be intersected with an arbitrary interval. In \eqref{eq:vaaler} we have $g=G=1$.

\begin{theorem}\label{thm:rec}
For every $\varepsilon>0$, uniformly over this family,
\[
|S|\lesssim_\varepsilon HT^{\beta+\varepsilon}\qquad(1\le H\le MT^{-\beta}).
\]
\end{theorem}

Since $|S|\le 4HM$ trivially, the theorem has content only for $M>T^{\beta}$, which is the case in \eqref{eq:vaaler}. It is proved in Sections~\ref{subsec:classical}--\ref{subsec:large}. Together with \eqref{eq:vaaler}, it bounds each sum \eqref{eq:psiblock} by $O_\varepsilon(T^{\beta+\varepsilon})$, and Theorem~\ref{thm:main} then follows from \eqref{eq:Epsi} and \eqref{eq:Dpsi}. The details are given in Section~\ref{subsec:proofmain}.

\emph{Parameters.} For the sums \eqref{eq:S} we put, as in \cite{CGZZ},
\begin{equation}\label{eq:params}
\rho=\frac{HT}{M},\qquad U=\frac{HT}{M^2},\qquad p=\log_T\rho,\qquad m=\log_TM,\qquad \tau=1-p,
\end{equation}
so that $H=T^{m-\tau}$, $U=T^{p-m}$, $\rho^2=THU$ and $U/H=T/M^2\ge 1$. Apart from the term $hT\kappa/M$, which does not depend on $n$, the phase in \eqref{eq:S} has size $\asymp\rho$, and its derivative with respect to $n$ has size $\asymp U$. The first step in Section~\ref{subsec:sieveproof} is the van der Corput $B$-process (see \cite[Chapter~3]{GK}), which replaces the sum over $n$ by a sum over the values $q$ of minus this derivative, so that $q\asymp U$. To find the new phase, put $\omega=M\alpha$. For $q>0$ the function $n\mapsto hT/(4(n+\omega))+qn$ has its stationary point at $n+\omega=(hT/(4q))^{1/2}$, and its value there is
\[
f_h(q)=\sqrt{Thq}-\omega q .
\]
This is the phase of the sums over $q$ in \cite[Section~3]{CGZZ}. Finally, the exponents $p$, $m$ and $\tau$ describe $\rho$, $M$ and $H$ on a logarithmic scale; the condition $H\le MT^{-\beta}$ in Theorem~\ref{thm:rec} means that $\tau\ge\beta$. The number $s$ is a threshold for $H$: if $H\le T^{s}$, a classical bound already gives Theorem~\ref{thm:rec}, because $(HT)^{2/7}\le T^{2(1+s)/7}=T^{\beta}$; see \eqref{eq:easy}.

\subsection{Classical bounds and the two main ranges}\label{subsec:classical}

In the proof of Theorem~\ref{thm:rec}, the letters $\rho$, $U$, $p$, $m$ and $\tau$ always denote the parameters \eqref{eq:params}. By \cite[Proposition~2.2]{CGZZ}, which holds for all sums \eqref{eq:S}, uniformly on subintervals,
\begin{equation}\label{eq:classical}
\frac{|S|}{H}\lesssim(HT)^{2/7}+\frac{M^2}{HT},\qquad
\frac{|S|}{H}\lesssim\rho^{13/31}M^{3/31}+\frac{M^2}{HT}.
\end{equation}
These bounds are obtained by applying the exponent pairs $(\frac27,\frac47)$ and $(\frac{13}{31},\frac{16}{31})$ to the sum over $n$; see \cite[Chapter~3]{GK}. The exponent pair $(\frac29,\frac{11}{18})=BA(\frac27,\frac47)$ gives in the same way
\begin{equation}\label{eq:classical3}
\frac{|S|}{H}\lesssim\rho^{2/9}M^{7/18}+\frac{M^2}{HT}.
\end{equation}
Note that $M^2/(HT)\le 1$. If $H\le T^{s}$, then $(HT)^{2/7}\le T^{2(1+s)/7}=T^{\beta}$. If $p\le 42/67$, then $\rho^{13/31}M^{3/31}\le T^{1293/4154}$, and $1293/4154=0.31126\ldots<\beta$. Hence, by \eqref{eq:classical},
\begin{equation}\label{eq:easy}
\frac{|S|}{H}\lesssim T^{\beta}\qquad\text{if}\quad H\le T^{s}\quad\text{or}\quad p\le\frac{42}{67}.
\end{equation}
The remaining cases are covered by the following two propositions. They are the analogues of \cite[Propositions~2.3 and~2.4]{CGZZ}, where the middle range is $T^{42/67}\le\rho\le T^{43/67}$.

\begin{proposition}[The middle range]\label{prop:middle}
Assume
\begin{equation}\label{eq:middle}
T^{42/67}\le\rho\le T^{61/94},\qquad H\ge T^{s},\qquad M\le T^{1/2}.
\end{equation}
Then, uniformly in the shifts, weights, and interval endpoints,
\[
\frac{|S|}{H}\lesssim_\varepsilon T^{\beta+\varepsilon}.
\]
\end{proposition}

\begin{proposition}[The larger-$\rho$ range]\label{prop:large}
Assume
\begin{equation}\label{eq:large}
\beta\le\tau\le\frac{33}{94},\qquad \tau+s\le m\le\frac12 .
\end{equation}
Then, uniformly in the shifts, weights, and interval endpoints,
\[
\frac{|S|}{H}\lesssim_\varepsilon T^{\beta+\varepsilon}.
\]
\end{proposition}

\begin{proof}[Proof of Theorem~\ref{thm:rec}]
If $H\le T^{s}$ or $p\le 42/67$, the conclusion follows from \eqref{eq:easy}. Otherwise $m>\tau+s$ and $\tau<25/67$, and $\tau\ge\beta$ because $H\le MT^{-\beta}$. For $33/94\le\tau<25/67$ the hypothesis \eqref{eq:middle} holds, and we apply Proposition~\ref{prop:middle}. For $\beta\le\tau\le 33/94$ we apply Proposition~\ref{prop:large}.
\end{proof}

\subsection{The middle range}\label{subsec:middle}

Throughout Sections~\ref{subsec:middle}--\ref{subsec:close} we assume \eqref{eq:middle}. The length of the shifts in Section~\ref{subsec:sieveproof} is
\begin{equation}\label{eq:J}
J=\Bigl\lceil\frac{(T\rho^{11})^{1/15}}{M}\Bigr\rceil .
\end{equation}
For $\rho=T^{43/67}$ this is the choice $J=\lceil T^{36/67}/M\rceil$ of \cite[(3.1)]{CGZZ}. As in \cite[(3.2), (3.3)]{CGZZ}, for $i\in\{1,\dots,6\}$ and $h,q>0$ let
\[
A_i(h,q)=\frac{f_h^{(i)}(q)}{i!}=\binom{1/2}{i}\sqrt{Th}\,q^{1/2-i}-\omega\,\mathbf{1}_{\{i=1\}}
\]
be the Taylor coefficients of the function $f_h$ from Section~\ref{subsec:reduction}, let $C_0\ge 1$ be a sufficiently large constant, let
\[
\Omega=\bigl([H,2H)\cap\mathbb{Z}\bigr)\times\bigl([U/100,2U]\cap\mathbb{Z}\bigr),
\]
and let
\begin{equation}\label{eq:B}
\begin{aligned}
B=\#\bigl\{\bigl((h,q),(h',q')\bigr)\in\Omega^2:\ &\|A_4'-A_4\|\le C_0J^{-4},\\
&|A_i'-A_i|\le C_0J^{-i}\ \ (i=5,6)\bigr\},
\end{aligned}
\end{equation}
where $A_i=A_i(h,q)$ and $A_i'=A_i(h',q')$. Thus $B$ counts the pairs of points of $\Omega$ whose Taylor coefficients of degrees $4$, $5$ and $6$ are close to each other, the coefficient of degree $4$ modulo~$1$.

\begin{lemma}\label{lem:csieve}
For the choice \eqref{eq:J},
\begin{equation}\label{eq:csieve}
\begin{aligned}
\frac{|S|}{H}\lesssim_\varepsilon T^{\varepsilon}\Bigl[&\bigl(T\rho^{4}M^{5}(MJ)^{6}\bigr)^{1/30}\Bigl(\frac{B}{HU}\Bigr)^{1/30}\\
&+\bigl(T\rho^{3}(MJ)^{6}\bigr)^{1/20}\Bigl(\frac{B}{HU}\Bigr)^{1/20}+\frac{MJ}{\sqrt{\rho}}+1\Bigr].
\end{aligned}
\end{equation}
\end{lemma}

In \cite[Lemma~3.1]{CGZZ} the bound is $T^{\varepsilon}[(T\rho^{8}(MJ)^{6})^{1/30}(B/(HU))^{1/30}+MJ/\sqrt{\rho}+1]$. Since $M^5/\rho^4=\rho/U^5$, the first term in \eqref{eq:csieve} is smaller than the first term there by the factor $(\rho/U^5)^{1/30}$.

\begin{proposition}\label{prop:close}
For the choice \eqref{eq:J},
\begin{equation}\label{eq:close}
B\lesssim_\varepsilon HU\,T^{\varepsilon}.
\end{equation}
\end{proposition}

Lemma~\ref{lem:csieve} and Proposition~\ref{prop:close} are proved in Sections~\ref{subsec:sieveproof} and~\ref{subsec:close}. Both proofs use the following relations between the parameters, which take the place of \cite[(3.11)]{CGZZ}; see also Table~\ref{tab:close}:
\begin{equation}\label{eq:sizes}
\begin{gathered}
J\ge T^{53/2010},\qquad \frac JU=\frac{MJ}{\rho}\lesssim T^{-101/1005},\\
\rho\Bigl(\frac JU\Bigr)^{7}=\frac{(MJ)^7}{\rho^6}\lesssim T^{-77/1005},\qquad \frac{\rho J^6}{U^6}=\frac{(MJ)^6}{\rho^5}\ge T^{1/94},\\
\frac{\rho}{U^5}=\frac{M^5}{\rho^4}\le T^{-1/134},\qquad \frac{\rho J^4}{U^5}=\frac{M(MJ)^4}{\rho^4}\ge T^{31/1880}.
\end{gathered}
\end{equation}
To verify them, note first that $(T\rho^{11})^{1/15}/M\ge T^{53/2010}$, because $M\le T^{1/2}$ and $\rho\ge T^{42/67}$. Hence $J\ge T^{53/2010}$ and
\begin{equation}\label{eq:MJ}
(T\rho^{11})^{1/15}\le MJ\le 2(T\rho^{11})^{1/15}.
\end{equation}
Inserting \eqref{eq:MJ} we find that
\[
\frac{MJ}{\rho}\lesssim\frac{T^{1/15}}{\rho^{4/15}},\qquad \frac{(MJ)^7}{\rho^6}\lesssim\frac{T^{7/15}}{\rho^{13/15}},\qquad \frac{(MJ)^6}{\rho^5}\ge\frac{T^{2/5}}{\rho^{3/5}}.
\]
The second and the third relation in \eqref{eq:sizes} now follow from $\rho\ge T^{42/67}$, and the fourth from $\rho\le T^{61/94}$. The fifth follows from $M\le T^{1/2}$ and $\rho\ge T^{42/67}$. Finally, $H\ge T^{s}$ means that $M\ge T^{1+s}/\rho$, so that $M(MJ)^4/\rho^4\ge T^{19/15+s}\rho^{-31/15}\ge T^{31/1880}$, by $\rho\le T^{61/94}$.

\begin{proof}[Proof of Proposition~\ref{prop:middle}]
By \eqref{eq:MJ}, $\rho\le T^{61/94}$ and $M\le T^{1/2}$,
\[
\begin{aligned}
\bigl(T\rho^{4}M^{5}(MJ)^{6}\bigr)^{1/30}&\lesssim\bigl(T^{7/5}\rho^{42/5}M^{5}\bigr)^{1/30}\le T^{(\frac75+\frac{42}{5}\cdot\frac{61}{94}+\frac52)/30}=T^{293/940},\\
\bigl(T\rho^{3}(MJ)^{6}\bigr)^{1/20}&\lesssim\bigl(T^{7/5}\rho^{37/5}\bigr)^{1/20}\le T^{(\frac75+\frac{37}{5}\cdot\frac{61}{94})/20}=T^{583/1880},\\
\frac{MJ}{\sqrt{\rho}}&\le 2T^{1/15}\rho^{7/30}\le 2T^{41/188}.
\end{aligned}
\]
The first exponent is $\beta$, and the other two are smaller. Hence the proposition follows from \eqref{eq:csieve} and \eqref{eq:close}, applied with smaller positive losses.
\end{proof}

\subsection{Proof of Lemma~\ref{lem:csieve}}\label{subsec:sieveproof}

We follow \cite[Section~3.2]{CGZZ}. With $f_h(q)=\sqrt{Thq}-\omega q$ as in Section~\ref{subsec:reduction}, let
\[
\mathcal{S}=\sup\Bigl|\sum_{h\approx H}c_h\sum_{q\in I_h\cap\mathbb{Z}}e\bigl(f_h(q)\bigr)\Bigr|,
\]
where the supremum is over all $|c_h|\le 1$ and all intervals $I_h\subseteq[U/100,2U]$. Abel summation and the van der Corput transformation of the sum over $n$ give \cite[(3.8)]{CGZZ}
\begin{equation}\label{eq:vdc}
\frac{|S|}{H}\lesssim\frac{M}{H\sqrt{\rho}}\,\mathcal{S}+\log(2+U)+\frac{M}{\sqrt{\rho}}.
\end{equation}
This is the $B$-process described in Section~\ref{subsec:reduction}. Its derivation in \cite{CGZZ} uses only the sizes of the derivatives of the phase in \eqref{eq:S}, which are the same for all the sums \eqref{eq:S} because $|\alpha|\le\frac14$. In particular, it depends neither on the choice of $J$ nor on the range of $\rho$.

Next, the sums over $q$ are expanded around short shifts. On the interval $[U/100,2U]$ we have $|f_h^{(7)}(q)|\asymp\sqrt{Th}\,q^{-13/2}\asymp\rho/U^7$, and $\rho J^7/U^7\lesssim 1$ by the third relation in \eqref{eq:sizes}. Hence Lemma~\ref{lem:shift} can be applied with $d=6$ to the functions $f_h$ on the intervals $I_h$, with the variable $q$ in place of $n$. The pairs $(h,q)$ that occur in the lemma belong to $\Omega$, and their Taylor coefficients are the numbers $A_i(h,q)$. Since there are at most $2H$ values of $h$, the lemma gives
\begin{equation}\label{eq:bilinear}
\mathcal{S}\lesssim HJ+\frac1J\,\sup\Bigl|\sum_{(h,q)\in\Omega}a_{h,q}\sum_{r=1}^{j}e\Bigl(\sum_{i=1}^{6}A_i(h,q)\,r^i\Bigr)\Bigr|,
\end{equation}
where the supremum is over $1\le j\le J$ and all $|a_{h,q}|\le 1$. Up to this point the argument is that of \cite{CGZZ}; it uses the choice of $J$ only through the third relation in \eqref{eq:sizes}.

The inner sum in \eqref{eq:bilinear} is a Weyl sum of degree six. We apply Lemma~\ref{lem:ls} with $d=5$ to the points $(A_1(h,q),\dots,A_6(h,q))$ with $(h,q)\in\Omega$, and to the coefficients $b_r=\mathbf{1}_{\{r\le j\}}$. The number of points is $L=\#\Omega\asymp HU$. Since $\sqrt{Th}\,q^{1/2}\asymp\rho$ on $\Omega$, there is a constant $C\ge 1$ such that $|A_5(h,q)|\le C\rho/U^5$ and $|A_6(h,q)|\le C\rho/U^6$, and we take
\[
D=\frac{C\rho}{U^5},\qquad D'=\frac{C\rho}{U^6}.
\]
The restriction of the fifth coordinate to the short interval $[-D,D]$ is the new element in comparison with \cite{CGZZ}. By the last two relations in \eqref{eq:sizes}, $J^{-4}\le D\le\frac14$ if $T$ is large, and by the fourth relation the number $\Delta=1+D'J^6$ of Lemma~\ref{lem:ls} satisfies $\Delta\asymp\rho J^6/U^6$. Every pair of points that is counted in Lemma~\ref{lem:ls} is counted in \eqref{eq:B}, so that the number of these pairs is at most $B$. Therefore \eqref{eq:ls} bounds the supremum in \eqref{eq:bilinear} by
\[
C_\varepsilon J^{1+\varepsilon}\,HU\Bigl[\Bigl(\frac{\Delta D}{HU}\cdot\frac{B}{HU}\Bigr)^{1/30}+\Bigl(\frac{\Delta}{HU}\cdot\frac{B}{HU}\Bigr)^{1/20}\Bigr].
\]
We insert this into \eqref{eq:bilinear} and \eqref{eq:vdc}. Since $U=\rho/M$ and $\rho^2=THU$,
\[
\Bigl[\frac{MU}{\sqrt{\rho}}\Bigr]^{30}\frac{\Delta D}{HU}\asymp T\rho^{4}M^{5}(MJ)^{6},\qquad
\Bigl[\frac{MU}{\sqrt{\rho}}\Bigr]^{20}\frac{\Delta}{HU}\asymp T\rho^{3}(MJ)^{6}.
\]
This gives the first two terms in \eqref{eq:csieve}. The term $HJ$ in \eqref{eq:bilinear} contributes $MJ/\sqrt{\rho}$, and the last two terms in \eqref{eq:vdc} are $O(T^{\varepsilon}+MJ/\sqrt{\rho})$. This proves Lemma~\ref{lem:csieve}.

\subsection{Proof of Proposition~\ref{prop:close}}\label{subsec:close}

For $\rho\le T^{43/67}$, $J=\lceil T^{36/67}/M\rceil$ and $123/1340$ in place of $s$, the bound \eqref{eq:close} is \cite[Proposition~3.1]{CGZZ}. The numbers $A_i(h,q)$, the set $\Omega$, the constant $C_0$ and the number $B$ in \eqref{eq:B} are defined as in \cite[(3.2), (3.3)]{CGZZ}; only the range of the parameters and the choice of $J$ are different. We keep the notation of \cite[Section~4]{CGZZ}. Thus the letters $a$, $b$, $k$, $\ell$, $n$ and $r$ have a new meaning in this subsection; in particular, $k$ and $\ell$ do not refer to exponent pairs. As in \cite{CGZZ} let
\[
a=\frac{U^6}{\rho J^6}=\frac{\rho^5}{(MJ)^6},\qquad b=\frac{U^5}{\rho J^5}=\frac{\rho^4}{(MJ)^5},
\]
so that $a=Ub/J$.

\emph{What the proof in \cite{CGZZ} uses.} The proof of \cite[Proposition~3.1]{CGZZ} in \cite[Section~4]{CGZZ} uses the parameters in the following ways only: through identities between them, such as $\rho^2=THU$, $U=\rho/M$ and $a=Ub/J$; through the inequality $H\le U$, which follows from $M\le T^{1/2}$; through the condition $C_0J^{-4}<\frac12$; and through the inequalities \cite[(4.2)]{CGZZ}, which are derived there from \cite[(2.3), (3.11)]{CGZZ}. The other ingredients of the proof, namely the divisor bound, a Fej\'er kernel, the van der Corput transformation and the estimates of Robert and Sargos~\cite{RS}, hold without any restriction on the parameters. In our situation, the inequalities \cite[(4.2)]{CGZZ} are replaced by
\begin{equation}\label{eq:sizes2}
\begin{gathered}
a\le T^{-1/94},\qquad UJ^{-5}=\frac{\rho b}{U^4}\le T^{-1/201},\qquad Jb\le T^{-31/1880},\\
HJb^2=\frac{\rho^9}{T(MJ)^9}\le T^{-2/47},\qquad \frac{H\rho b^3}{U^3}=\frac{\rho^{11}M^4}{T(MJ)^{15}}\le 1,\\
\frac{\rho b^2}{U^3}=\frac{\rho^6M^3}{(MJ)^{10}}\le T^{-1/402},\qquad \frac{\rho b}{HU^3}=\frac{TM^2\rho}{(MJ)^5}\le T^{-1/201}.
\end{gathered}
\end{equation}
They are verified in the same way as \eqref{eq:sizes}. The first and the third inequality are the fourth and the sixth relation in \eqref{eq:sizes}. For the second, fourth, sixth and seventh we insert the lower bound for $MJ$ in \eqref{eq:MJ}, which gives
\[
\begin{gathered}
\frac{\rho M^4}{(MJ)^5}\le\frac{M^4}{T^{1/3}\rho^{8/3}},\qquad
\frac{\rho^9}{T(MJ)^9}\le\frac{\rho^{12/5}}{T^{8/5}},\\
\frac{\rho^6M^3}{(MJ)^{10}}\le\frac{M^3}{T^{2/3}\rho^{4/3}},\qquad
\frac{TM^2\rho}{(MJ)^5}\le\frac{T^{2/3}M^2}{\rho^{8/3}},
\end{gathered}
\]
and then we use $M\le T^{1/2}$ and $T^{42/67}\le\rho\le T^{61/94}$. The fifth inequality holds because $(MJ)^{15}\ge T\rho^{11}$ and $M^4\le T^2$. It is the critical one, and it is the reason for the choice \eqref{eq:J}.

Table~\ref{tab:close} compares \eqref{eq:sizes} and \eqref{eq:sizes2} with the corresponding bounds in \cite[(3.11), (4.2)]{CGZZ}, and it indicates where each of them is used. Every bound that is a negative power of $T$ in \cite{CGZZ} is also a negative power of $T$ here, and the bound $1$ for $H\rho b^3/U^3$ is the same in both papers. Moreover, $C_0J^{-4}<\frac12$ if $T$ is large, because $J\ge T^{53/2010}$. Therefore the argument of \cite[Section~4]{CGZZ} applies to our parameters without change; only the implied constants and the lower bound for $T$ may be different.

\begin{table}[ht]
\caption{Upper bounds for the quantities in the first column under the assumptions of \cite[Section~3]{CGZZ} and under \eqref{eq:middle}, \eqref{eq:J}. The first part corresponds to \eqref{eq:sizes}, the second to \eqref{eq:sizes2}; note that $a$ and $Jb$ are the reciprocals of $\rho J^6/U^6$ and $\rho J^4/U^5$. The fourth column gives the parameters at which our bound is attained, with $\rho_0=T^{42/67}$ and $\rho_1=T^{61/94}$, and the last column the places where the bound is used.}\label{tab:close}
\small
\renewcommand{\arraystretch}{1.2}
\begin{tabular}{@{}lllll@{}}
\hline
quantity & \cite{CGZZ} & here & worst case & used in\\
\hline
$J^{-1}$ & $T^{-5/134}$ & $T^{-53/2010}$ & $\rho_0$, $M=T^{1/2}$ & $C_0J^{-4}<\frac12$\\
$J/U$ & $T^{-6/67}$ & $T^{-101/1005}$ & $\rho_0$ & \cite[(4.5)]{CGZZ}\\
$\rho J^7/U^7$ & $1$ & $T^{-77/1005}$ & $\rho_0$ & Lemma~\ref{lem:shift} in Sect.~\ref{subsec:sieveproof}\\
$\rho/U^5$ & $T^{-1/134}$ & $T^{-1/134}$ & $\rho_0$, $M=T^{1/2}$ & $D\le\frac14$ in Sect.~\ref{subsec:sieveproof}\\
\hline
$a$ & $T^{-1/67}$ & $T^{-1/94}$ & $\rho_1$ & \cite[(4.4)--(4.6)]{CGZZ}; Sect.~\ref{subsec:sieveproof}\\
$UJ^{-5}$ & $T^{-3/67}$ & $T^{-1/201}$ & $\rho_0$, $M=T^{1/2}$ & \cite[(4.4), (4.6)]{CGZZ}\\
$Jb$ & $T^{-43/1340}$ & $T^{-31/1880}$ & $\rho_1$, $H=T^{s}$ & \cite[(4.5)]{CGZZ}; Sect.~\ref{subsec:sieveproof}\\
$HJb^2$ & $T^{-4/67}$ & $T^{-2/47}$ & $\rho_1$ & \cite[(4.5)]{CGZZ}\\
$H\rho b^3/U^3$ & $1$ & $1$ & $M=T^{1/2}$ & \cite[(4.11)]{CGZZ}\\
$\rho b^2/U^3$ & $T^{-3/134}$ & $T^{-1/402}$ & $\rho_0$, $M=T^{1/2}$ & \cite[(4.18)]{CGZZ}\\
$\rho b/(HU^3)$ & $T^{-3/67}$ & $T^{-1/201}$ & $\rho_0$, $M=T^{1/2}$ & \cite[(4.11), (4.18)]{CGZZ}\\
\hline
\end{tabular}
\end{table}

In outline, the conditions on $A_5$ and $A_6$ confine $(h',q')$ to a small neighborhood of a curve through $(h,q)$, and the condition on $A_4$ produces an integer $\ell$. The pairs with $\ell=0$ are counted with the divisor bound. For $\ell\ne 0$ the count becomes a count of integer points near curves, which is carried out with exponential sums. We now go through the argument of \cite[Section~4]{CGZZ} in order to indicate where the inequalities \eqref{eq:sizes2} are used. Let $T$ be large. Consider a pair of points $(h,q)$, $(h',q')$ that is counted in \eqref{eq:B}, with $A_i$ and $A_i'$ as there, and write $k=q'-q$ and $w=h'-h$. Since $C_0J^{-4}<\frac12$, there is a unique integer $\ell$ with $|A_4'-A_4-\ell|\le C_0J^{-4}$. By the relations $A_5=-\frac43qA_6$ and $A_4=-\frac{10}{7}qA_5$, the conditions on $A_5'-A_5$ and $A_6'-A_6$ give \cite[(4.4)]{CGZZ}
\[
|k|\lesssim Ua,\qquad h+w=h\Bigl(1+\frac kq\Bigr)^{9}+O(Hb),\qquad A_4'-A_4=\frac kq\,A_4+O(UJ^{-5}).
\]

\emph{The case $\ell=0$.} Here the conditions on $A_4'-A_4$ and $A_5'-A_5$ give $|k|\lesssim Ub$, $|w|\lesssim Hb$ and $|qw-7hk|\lesssim HJb+HUb^2$. Since $H\le U$, $Jb=o(1)$ and $Ub^2=aJb$, a pair with $kw=0$ is diagonal. For the other pairs one chooses $h$, $k$ and the integer $qw-7hk$, and then the divisor bound leaves $O_\varepsilon(T^{\varepsilon})$ possibilities for $q$ and $w$. Hence the number $B^0$ of pairs with $\ell=0$ satisfies \cite[(4.5)]{CGZZ}
\[
\frac{B^0}{HU}\lesssim_\varepsilon T^{\varepsilon}\bigl(1+b+HJb^2+HUb^3\bigr)\lesssim_\varepsilon T^{\varepsilon},
\]
because $HUb^3=aHJb^2$.

\emph{The case $\ell\ne 0$.} We may assume that $k>0$; then $\ell<0$, because $A_4<0$. We consider the pairs with $k\approx Ur$, for $O(\log T)$ dyadic values of $r$, and as in \cite{CGZZ} we put $n=-\ell$, $K=Ur$ and $Z=\rho r/U^4$. Then $n\asymp Z$, $w\asymp Hr$, $U^4/\rho\lesssim r\lesssim a$ and $b/r=o(1)$, as in \cite[(4.6)]{CGZZ}; for the last relation note that $b/r\lesssim\rho b/U^4=UJ^{-5}$. The conditions are then inverted in \cite[(4.7)--(4.9)]{CGZZ}. With
\[
f_{k,n}(h)=\bigl(h^{1/9}+\mathcal{A}_{k,n}^{1/9}\bigr)^{9}-h,\qquad\mathcal{A}_{k,n}=\frac{16384\,n^2k^7}{25\,T},\qquad D_r=\frac{Hb^2}{r},
\]
they imply that $|w-f_{k,n}(h)|\lesssim D_r$, and they leave $O(Ub/r)$ possibilities for $q$ when $h$, $w$, $k$ and $n$ are fixed; note that $Ub/r\gtrsim J$. Let $B_r$ be the number of these pairs, and let $Q=\lfloor c\rho br/U^3\rfloor$ with the small constant $c$ of \cite[(4.10)]{CGZZ}. If $Q<2$, the trivial count of $O(1+D_r)$ values of $w$ gives
\[
\frac{B_r}{HU}\lesssim\frac{\rho br}{U^3}+\frac{H\rho b^3}{U^3}\lesssim 1 .
\]
If $Q\ge 2$, the condition on $w$ is detected by a Fej\'er kernel of length $Q$. This requires \cite[(4.11)]{CGZZ}
\[
QD_r\le\frac{cH\rho b^3}{U^3}\le c,\qquad \frac{Q}{Hr}\lesssim\frac{\rho b}{HU^3}\le T^{-1/201}.
\]
The nonzero frequencies $v$ of the kernel are divided into the ranges $V\le v<\min(2V,Q)$ with $V=2^j<Q$. By the last display, $V\lesssim H$, $rV\le H$ and $Y=HrV\ge Hr\gtrsim T^{1/201}$. For each $v$ the sum over $h$ is transformed by the van der Corput lemma, and the resulting sums are estimated by a double large sieve with the two spacing counts $B_x\lesssim_\varepsilon T^{\varepsilon}(ZK+(ZK)^2/Y)$ and $B_y\lesssim_\varepsilon T^{\varepsilon}V^2$, which come from the work of Robert and Sargos~\cite{RS}; see \cite[(4.12)--(4.17)]{CGZZ}. The result is \cite[(4.18)]{CGZZ}
\[
\frac{B_r}{HU}\lesssim_\varepsilon T^{\varepsilon}\Bigl(1+\frac{Q}{\sqrt{ZK}}+\sqrt{\frac{Q}{Hr}}\Bigr)\lesssim_\varepsilon T^{\varepsilon}\Bigl(1+\sqrt{\frac{\rho b^2}{U^3}}+\sqrt{\frac{\rho b}{HU^3}}\Bigr)\lesssim_\varepsilon T^{\varepsilon}.
\]
Summing over $r$ we obtain \eqref{eq:close}. Every step uses only the inequalities \eqref{eq:sizes} and \eqref{eq:sizes2} and $H\le U$. In Section~\ref{subsec:sharp} we show that the choice \eqref{eq:J} is optimal for this argument.

\subsection{The larger-\texorpdfstring{$\rho$}{rho} range}\label{subsec:large}

In this subsection we assume \eqref{eq:large} and prove Proposition~\ref{prop:large}. The letters $N$, $R$ and $W$ now denote the parameters of \cite[Section~5]{CGZZ}. As in \cite[(5.1), (5.3)]{CGZZ} let
\[
N_{\mathrm{I}}=2^{18/35}M^{57/35}H^{-22/35}T^{-17/35},\quad\
N_{\mathrm{II}}=\Bigl(\frac{M^6}{TH^5}\Bigr)^{1/7},\quad\
N_{\mathrm{III}}=2\Bigl(\frac{M^6}{T^2H}\Bigr)^{1/3},
\]
\[
N=\begin{cases}\lceil N_{\mathrm{I}}\rceil,& HM^9\ge T^4,\\[2pt] \lceil\min(N_{\mathrm{II}},N_{\mathrm{III}})\rceil,& HM^9<T^4,\end{cases}
\qquad R=\Bigl\lceil\sqrt{\frac{C_RM^3}{2NT}}\Bigr\rceil,\qquad W=\frac HR,
\]
where $C_R$ is the constant of \cite[(5.3)]{CGZZ}. Here $N$ is the length of the short intervals into which the range of $n$ is cut, and $R$ controls the denominators of the rational approximations that are attached to these intervals; we refer to \cite[Section~5]{CGZZ} for the details. By \cite[Lemma~5.1]{CGZZ},
\begin{equation}\label{eq:watt}
\frac{|S|}{H}\lesssim_\varepsilon T^{\varepsilon}\Bigl[\frac{M}{\sqrt{HN}}\,W^{2/9}+\frac{M}{N^{3/4}H^{1/4}}\,W^{1/3}\Bigr],
\end{equation}
provided that $\frac{209}{670}\le\tau\le\frac{24}{67}$ and $\tau+\frac{123}{1340}\le m\le\frac12$. The range \eqref{eq:large} is not contained in this one, because $\beta<209/670$ and $s<123/1340$.

\emph{What the proof in \cite{CGZZ} uses.} In the proof of \eqref{eq:watt} in \cite[Section~5]{CGZZ}, Lemmas~8.1 and~8.2 of Watt~\cite{Watt} are applied with the parameters $(M_*,H_*,T_*)=(2M,2H,2T)$, with $N$ and $R$ as above, and with the functions $f$ and $\Phi$ of \cite[(5.1)]{CGZZ}. The restrictions on $\tau$ and $m$ are used only to verify the hypotheses of these lemmas and the bound \cite[(5.16)]{CGZZ} for an error term; this is done in \cite[Sections~5.1 and~5.2]{CGZZ}. The rest of the proof uses only identities between the parameters. Some of the hypotheses do not depend on the range of $\tau$ and $m$. The condition \cite[(2.6)]{Watt}, which reads $2(R-1)^2NT_*<C_3M_*^3\le 2R^2NT_*$, holds by the definition of $R$. The inequality $T_*/M_*^2\ge\frac12$ follows from $M\le T^{1/2}$, and $M_*\ge 4$ holds because $m\ge\tau+s>0$ and $T$ is large. The properties of $\Phi$ in \cite[(5.2)]{CGZZ} depend only on $|\alpha|\le\frac14$ and $M\le T^{1/2}$. The remaining hypotheses are the inequalities in Table~\ref{tab:watt}, which we now verify under the assumption \eqref{eq:large}.

Since $H=T^{m-\tau}$, the numbers $N_{\mathrm{I}}$, $N_{\mathrm{II}}$ and $N_{\mathrm{III}}$ are, up to constant factors, the powers of $T$ with the exponents $m+\frac{22\tau-17}{35}$, $\frac{m+5\tau-1}{7}$ and $\frac{5m+\tau-2}{3}$. By $m\ge\tau+s$ and $\tau\ge\beta$, each of these exponents is at least $13/120>1/10$, so that taking ceilings changes $N$ only by a fixed factor. The condition $HM^9\ge T^4$ is the same as $10m-\tau\ge 4$. The other inequalities in Table~\ref{tab:watt} read
\begin{equation}\label{eq:sizes3}
\begin{aligned}
\frac{HNT}{M^3}&\lesssim T^{(8-33\tau)/70}+T^{(1-3\tau-s)/3}\lesssim T^{-49/5640},\\
\frac{H}{N}&\lesssim T^{(17-57\tau)/35}+T^{(2-6\tau-2s)/3}\lesssim T^{-49/2820},\\
\frac{R}{H}&\lesssim T^{3(8\tau-3)/35}+T^{(3\tau-2s-1)/6}\lesssim T^{-27/1645},\\
\frac{N}{M}&\lesssim T^{(22\tau-17)/35}+T^{(-\tau-6s-1)/7}\lesssim T^{-436/1645},\\
\frac{HN^3T}{M^4}&\lesssim T^{(31\tau-16)/35}+T^{(\tau-1)/5}\lesssim T^{-61/470},\\
\frac{H^2N^3}{MR^4}&\lesssim T^{(8\tau-3)/7}+T^{(8\tau-3)/16}\lesssim T^{-9/752}.
\end{aligned}
\end{equation}
In the third and the last line, $R$ is replaced by $(M^3/(NT))^{1/2}$; this is justified below. We explain the two terms in each line. If $HM^9\ge T^4$, then $N\asymp N_{\mathrm{I}}$ and $m\ge(4+\tau)/10$. Inserting $N\asymp N_{\mathrm{I}}$ gives the first term in each line: from the second line on, this is an identity, and in the first line we also use $m\ge(4+\tau)/10$. If $HM^9<T^4$, then $N\asymp\min(N_{\mathrm{II}},N_{\mathrm{III}})$ and $\tau+s\le m<(4+\tau)/10$. In the first, fourth, fifth and last line the quantity on the left increases with $N$, and it is bounded by the second term. Here we use $N\lesssim N_{\mathrm{III}}$ in the first and the fifth line, together with the lower bound for $m$ in the first line and the upper bound in the fifth; $N\lesssim N_{\mathrm{II}}$ and the lower bound for $m$ in the fourth line; and $\min(N_{\mathrm{II}},N_{\mathrm{III}})\le N_{\mathrm{II}}^{7/16}N_{\mathrm{III}}^{9/16}$ in the last line. In the second and the third line the quantity on the left decreases when $N$ increases. It is therefore bounded by the sum of its values for $N=N_{\mathrm{II}}$ and $N=N_{\mathrm{III}}$. The value $N_{\mathrm{II}}$ gives at most the first term, because of the upper bound for $m$, and the value $N_{\mathrm{III}}$ gives at most the second term, by the lower bound for $m$. The final bounds in \eqref{eq:sizes3} follow from $\beta\le\tau\le 33/94$, since all the exponents are monotonic in $\tau$; the value of $\tau$ at which each of them is attained is given in Table~\ref{tab:watt}.

By the first line of \eqref{eq:sizes3}, $M^3/(NT)\gtrsim HT^{49/5640}$. Hence $R\asymp(M^3/(NT))^{1/2}$, as was used above, and $R$ tends to infinity. Consequently $H/R^2\asymp HNT/M^3$ and $HN^2/(MR^2)\asymp HN^3T/M^4$, and $W=H/R$ tends to infinity by the third line. Thus all the hypotheses in Table~\ref{tab:watt} hold under the assumption \eqref{eq:large}, and the proof of \cite[Lemma~5.1]{CGZZ} gives \eqref{eq:watt} in this range.

\begin{table}[ht]
\caption{The hypotheses of Watt's lemmas that depend on $\tau$ and $m$, and the bound \cite[(5.16)]{CGZZ}. Here $(M_*,H_*,T_*)=(2M,2H,2T)$, and $C_3$ and $\varepsilon_0$ are constants from \cite{Watt}. The third and the fourth column give the bounds under the assumptions of \cite[Section~5]{CGZZ} and under \eqref{eq:large}, and the last column the value of $\tau$ at which our bound is attained.}\label{tab:watt}
\small
\renewcommand{\arraystretch}{1.2}
\setlength{\tabcolsep}{4.5pt}
\begin{tabular}{@{}lllll@{}}
\hline
quantity & needed for & \cite{CGZZ} & here & worst case\\
\hline
$N_{\mathrm{I}}^{-1}$, $N_{\mathrm{II}}^{-1}$, $N_{\mathrm{III}}^{-1}$ & size of $N$ & $T^{-443/4020}$ & $T^{-13/120}$ & $\tau=\beta$\\
$HNT/M^3$ & $R\asymp(M^3/(NT))^{1/2}$ & $T^{-37/4020}$ & $T^{-49/5640}$ & $\tau=\beta$\\
$H/R^2$ & $4C_3\sqrt{H_*}\le R$; \cite[(5.16)]{CGZZ} & $T^{-37/4020}$ & $T^{-49/5640}$ & $\tau=\beta$\\
$R/H$ & $R\le\varepsilon_0H_*$ & $T^{-27/2345}$ & $T^{-27/1645}$ & $\tau=33/94$\\
$H/N$ & \cite[(2.4), (2.7)--(2.9)]{Watt} & $T^{-37/2010}$ & $T^{-49/2820}$ & $\tau=\beta$\\
$N/M$ & \cite[(2.4), (2.7)--(2.9)]{Watt} & $T^{-611/2345}$ & $T^{-436/1645}$ & $\tau=33/94$\\
$HN^3T/M^4$ & \cite[(2.4), (2.7)--(2.9)]{Watt}; \cite[(5.16)]{CGZZ} & $T^{-43/335}$ & $T^{-61/470}$ & $\tau=33/94$\\
$H^2N^3/(MR^4)$ & \cite[(8.14)]{Watt} & $T^{-9/1072}$ & $T^{-9/752}$ & $\tau=33/94$\\
\hline
\end{tabular}
\end{table}

It remains to evaluate \eqref{eq:watt}. As in \cite[(5.5)]{CGZZ},
\[
\begin{aligned}
\frac{M}{\sqrt{HN}}\,W^{2/9}&\asymp T^{1/9}M^{2/3}H^{-5/18}N^{-7/18}\\
&\lesssim
\begin{cases}T^{(9+\tau)/30},& HM^9\ge T^4,\\[2pt] T^{(2m+1)/6}+T^{(10-7m+4\tau)/27},& HM^9<T^4,\end{cases}
\end{aligned}
\]
and the second term in \eqref{eq:watt} is equal to $(H/M)^{1/2}=T^{-\tau/2}$ times the power $\frac32$ of the first term. Hence, if the first term is $O(T^{\beta})$, then the second is $O(T^{(3\beta-\tau)/2})$, and this is $O(T^{\beta})$ because $\tau\ge\beta$.

If $HM^9\ge T^4$, then $(9+\tau)/30\le\beta$, because $\tau\le 33/94=30\beta-9$. If $HM^9<T^4$, then $10m-\tau<4$ and therefore $(2m+1)/6<(9+\tau)/30\le\beta$. Thus Proposition~\ref{prop:large} follows from \eqref{eq:watt}, except when $HM^9<T^4$ and $(10-7m+4\tau)/27>\beta$. In this case we use \eqref{eq:classical3} instead. The exponent of $\rho^{2/9}M^{7/18}$ is
\[
\frac{2(1-\tau)}{9}+\frac{7m}{18}=\frac79-\frac32\cdot\frac{10-7m+4\tau}{27}<\frac79-\frac32\,\beta<\beta,
\]
where the last inequality holds because $\beta>14/45$. This completes the proof of Proposition~\ref{prop:large}.

\subsection{Proof of Theorem~\ref{thm:main}}\label{subsec:proofmain}

The argument is that of \cite[Proposition~2.1]{CGZZ}, with Theorem~\ref{thm:rec} in place of \cite[Theorem~1.2]{CGZZ}. Let $T\ge 2$, $1\le M\le T^{1/2}$, $|\alpha|\le\frac14$ and $|\kappa|\le 1$, and let $I\subseteq[M,2M)$ be an interval. If $M\ge 2T^{\beta}$, Theorem~\ref{thm:rec} with $g=G=1$ bounds each supremum in \eqref{eq:vaaler} by $O_\varepsilon(HT^{\beta+\varepsilon})$, and there are $O(\log T)$ values of $H$. If $M<2T^{\beta}$, the sum \eqref{eq:psiblock} is trivially $O(T^{\beta})$. Hence, in both cases,
\begin{equation}\label{eq:psi}
\sum_{n\in I\cap\mathbb{Z}}\psi\Bigl(\frac TM\,F\Bigl(\frac nM\Bigr)\Bigr)\lesssim_\varepsilon T^{\beta+\varepsilon}.
\end{equation}
By \eqref{eq:Epsi} and the dyadic subdivision described in Section~\ref{subsec:reduction}, $E(X)$ is, up to $O(1)$, a combination with coefficients $\pm4$ of $O(\log X)$ sums of the form \eqref{eq:psiblock} with $T=X$. By \eqref{eq:psi}, each of them is $O_\varepsilon(X^{\beta+\varepsilon})$, and therefore $E(X)\lesssim_\varepsilon X^{\beta+\varepsilon}$.

The sums that occur in the divisor problem are the special case $\alpha=\kappa=0$ of those in the circle problem, so no further estimates are needed. By \eqref{eq:Dpsi} and a dyadic subdivision, $\Delta(X)$ is, up to $O(1)$, a combination with coefficients $-2$ of $O(\log X)$ sums of the form \eqref{eq:psiblock} with $T=4X$, $\alpha=\kappa=0$ and $1\le M\le\sqrt{X}$; for these parameters $\frac TM F(\frac nM)=\frac Xn$. By \eqref{eq:psi}, each of these sums is $O_\varepsilon(X^{\beta+\varepsilon})$. This gives $\Delta(X)\lesssim_\varepsilon X^{\beta+\varepsilon}$ and completes the proof of Theorem~\ref{thm:main}.

If \cite[Theorem~1.2]{CGZZ} is used in place of Theorem~\ref{thm:rec}, the same argument gives $\Delta(X)\lesssim_\varepsilon X^{209/670+\varepsilon}$, as was mentioned in the introduction.

\subsection{Optimality of the shift length}\label{subsec:sharp}

This subsection is not needed for the proof of Theorem~\ref{thm:main}. Its purpose is to show that the argument of Sections~\ref{subsec:middle}--\ref{subsec:close} is tight. First, the length \eqref{eq:J} of the shifts cannot be reduced: for shorter shifts, the number of close pairs exceeds the number of points by a positive power of $T$, so that Proposition~\ref{prop:close} fails. Second, at the boundary point $\rho=T^{61/94}$, $M=T^{1/2}$ of the middle range, where the exponent $293/940$ arises, no choice of the shift length makes the first term of the large sieve smaller than $T^{293/940}$ when the count \eqref{eq:B} is used. In this sense, the exponent of Theorem~\ref{thm:main} is the limit of the present argument; it is not a limit for the circle and divisor problems themselves. Let $B(J)$ denote the number \eqref{eq:B} for an arbitrary positive integer $J$, and recall that $L=\#\Omega\asymp HU$.

\emph{A lower bound for $B(J)$.} On $\Omega$ we have $|A_5|\lesssim\rho/U^5$ and $|A_6|\lesssim\rho/U^6$. We divide $\mathbb{T}$, which contains the values of $A_4$ modulo~$1$, into $J^4$ arcs of length $J^{-4}$, and the ranges of $A_5$ and $A_6$ into intervals of length at most $J^{-5}$ and $J^{-6}$, respectively. This gives $O(J^4(1+\rho J^5/U^5)(1+\rho J^6/U^6))$ cells, and two points of $\Omega$ in the same cell form a pair that is counted in $B(J)$. By the Cauchy--Schwarz inequality,
\begin{equation}\label{eq:packing}
\frac{B(J)}{L}\gtrsim\max\Bigl\{1,\ \frac{L}{J^4(1+\rho J^5/U^5)(1+\rho J^6/U^6)}\Bigr\}.
\end{equation}
If $\rho J^5/U^5\ge 1$ and $\rho J^6/U^6\ge 1$, the second entry is comparable to $HU^{12}/(\rho^2J^{15})=H\rho b^3/U^3$, where $b=U^5/(\rho J^5)$ as in Section~\ref{subsec:close}. Thus the critical quantity in \eqref{eq:sizes2} is, up to a constant factor, a lower bound for the average number $B(J)/L$ of points that are close to a given point. In particular, under these two conditions, \eqref{eq:close} can only hold if $H\rho b^3/U^3=O_\varepsilon(T^{\varepsilon})$.

\emph{An example.} Let $P$ be a large integer, and let
\[
T=P^{94},\qquad M=P^{47},\qquad H=U=P^{14},\qquad \rho=P^{61}.
\]
These parameters satisfy \eqref{eq:middle}, with the largest values $\rho=T^{61/94}$ and $M=T^{1/2}$ allowed there. Moreover, $L\asymp P^{28}$, $\rho/U^5=P^{-9}$ and $\rho/U^6=P^{-23}$, and \eqref{eq:J} gives $J_0=P^4$. Let $0\le\theta\le\frac16$ and $J=\lfloor P^{4-\theta}\rfloor$, and let $B^0(J)$ be the number of the pairs counted in $B(J)$ for which the integer $\ell$ of Section~\ref{subsec:close} vanishes. We claim that
\begin{equation}\label{eq:sharp}
LP^{15\theta}\lesssim B(J)-B^0(J)\le B(J)\lesssim_\varepsilon T^{\varepsilon}LP^{15\theta}.
\end{equation}
For $\theta=\frac16$ this gives $B(J)\gtrsim T^{5/188}L$, so that \eqref{eq:close} fails when the shifts are shorter than \eqref{eq:J} by the factor $P^{1/6}$.

For the lower bound we restrict both points to the square $[P^{14},2P^{14})^2$, which contains $P^{28}$ points of $\Omega$. Since $\binom{1/2}{5}=\frac{7}{256}$ and $\binom{1/2}{6}=-\frac{21}{1024}$, we have $0<A_5<P^{-9}/20$ and $-P^{-23}/20<A_6<0$ on this square. Hence at most $\frac14P^{28-15\theta}$ of the cells above contain points of the square, and the Cauchy--Schwarz inequality gives at least $3P^{28+15\theta}$ ordered pairs of distinct points of the square that are counted in $B(J)$. The pairs with $\ell=0$ are treated as in the case $\ell=0$ in Section~\ref{subsec:close}, with $k=q'-q$ and $w=h'-h$. Since $|A_4|\asymp P^5$ and $|A_5|\asymp P^{-9}$ on the square, the conditions on $A_4'-A_4$ and $A_5'-A_5$, together with $A_4=-\frac{10}{7}qA_5$ and the fact that $A_4$ is a constant multiple of $h^{1/2}q^{-7/2}$, give
\[
|k|\lesssim P^{3+5\theta},\qquad |qw-7hk|\lesssim P^{7+4\theta};
\]
here we used $6\theta\le 1$. If $kw=0$, these bounds imply $k=w=0$. Otherwise $h$, $k$ and $qw-7hk$ determine $(q,w)$ up to $O_\varepsilon(P^{\varepsilon})$ possibilities, by the divisor bound. Hence at most $O_\varepsilon(P^{24+9\theta+\varepsilon})$ pairs of distinct points have $\ell=0$, and the lower bound in \eqref{eq:sharp} follows. For the upper bound we sort the points of $\Omega$ into the cells of a grid with the widths $J_0^{-4}$, $J_0^{-5}$ and $J_0^{-6}$. A point of a given cell can form a pair counted in $B(J)$ only with points of $O(P^{15\theta})$ cells, and two points of the same cell form a pair counted in $B(J_0)$. Hence $B(J)\lesssim P^{15\theta}B(J_0)$, and \eqref{eq:close} gives the upper bound.

\emph{Stationary points.} For $J=J_0$ the lower bound in \eqref{eq:sharp} remains true if both points are restricted to the range that comes from the stationary points of the van der Corput transformation. Let $\alpha=0$, and write $x=h/H$ and $y=q/U$. The stationary point $(hT/(4q))^{1/2}$ lies in $(M,2M)$ if and only if $x/16<y<x/4$. The map
\[
(x,y)\longmapsto\bigl(P^{9}A_5,\,-P^{23}A_6\bigr)=\Bigl(\frac{7}{256}\sqrt{x}\,y^{-9/2},\ \frac{21}{1024}\sqrt{x}\,y^{-11/2}\Bigr)
\]
is injective, because the quotient of its components is a multiple of $y$. The absolute value of its Jacobian determinant is $\frac{147}{524288}\,y^{-11}$, which equals $\frac{147}{256}$ at $y=\frac12$. Let $\Gamma$ be a closed rectangle of positive area in $\{1<x<2,\ x/16<y<x/4\}$ on which the Jacobian is less than $\frac34$; there are such rectangles near $(2,\frac12)$. The image of $\Gamma$ has area at most $\frac34|\Gamma|$. Hence it meets at most $(\frac34|\Gamma|+o(1))P^{12}$ rectangles with the side lengths $P^{-11}$ and $P^{-1}$, which correspond to the widths $J_0^{-5}$ and $J_0^{-6}$. Together with the $P^{16}$ arcs for $A_4$, this gives fewer cells than the $(|\Gamma|+o(1))P^{28}$ points $(h,q)$ with $(x,y)\in\Gamma$. As before, we obtain $\gtrsim_{\Gamma}P^{28}$ pairs of such points with $\ell\ne 0$ that are counted in $B(J_0)$.

\emph{The large sieve.} Finally, at the parameters of the example no choice of $J$ makes the first term of the large sieve smaller. For every $J$ for which the proof of Lemma~\ref{lem:csieve} applies, the contribution of the first term in \eqref{eq:ls} is, up to a constant factor and the factor $T^{\varepsilon}$,
\[
\sqrt{\rho}\,\Bigl(\frac{\rho}{U^5}\cdot\frac{\Delta B(J)}{L^2}\Bigr)^{1/30},\qquad \Delta=1+\frac{\rho J^6}{U^6}.
\]
We have $\Delta B(J)/L\gtrsim P$ for every $J\ge 1$. Indeed, if $J\ge J_0$, then $B(J)\ge L$ and $\Delta\ge\rho J_0^6/U^6=P$. If $J\le J_0$, then \eqref{eq:packing} gives
\[
\frac{\Delta B(J)}{L}\gtrsim\frac{L}{J^4(1+\rho J^5/U^5)}\ge\frac{L}{J_0^4(1+\rho J_0^5/U^5)}\asymp P .
\]
Hence the contribution is $\gtrsim P^{61/2}(P^{-9}\cdot P\cdot P^{-28})^{1/30}=T^{293/940}$ for every $J$, while for $J=J_0$ it is $O_\varepsilon(T^{293/940+\varepsilon})$ by Proposition~\ref{prop:close}. Thus, with the count \eqref{eq:B}, the exponent $293/940$ is the best that the first term of the large sieve can give at this point, and the choice \eqref{eq:J} attains it.

\end{document}